\documentstyle[12pt,amssymb]{article}
\begin{document}

\baselineskip15.5pt

\newtheorem{definition}{Definition $\!\!$}[section]
\newtheorem{prop}[definition]{Proposition $\!\!$}
\newtheorem{lem}[definition]{Lemma $\!\!$}
\newtheorem{corollary}[definition]{Corollary $\!\!$}
\newtheorem{theorem}[definition]{Theorem $\!\!$}
\newtheorem{example}[definition]{\it Example $\!\!$}
\newtheorem{remark}[definition]{\it Remark $\!\!$}

\newcommand{\nc}[2]{\newcommand{#1}{#2}}
\newcommand{\rnc}[2]{\renewcommand{#1}{#2}}

\rnc{\theequation}{\thesection.\arabic{equation}}

\def\eps{{\epsilon}}
\nc{\qpb}{quantum principal bundle}
\nc{\bpr}{\begin{prop}}
\nc{\bth}{\begin{theorem}}
\nc{\ble}{\begin{lem}}
\nc{\bco}{\begin{corollary}}
\nc{\bre}{\begin{remark}}
\nc{\bex}{\begin{example}}
\nc{\bde}{\begin{definition}}
\nc{\ede}{\end{definition}}
\nc{\epr}{\end{prop}}
\nc{\ethe}{\end{theorem}}
\nc{\ele}{\end{lem}}
\nc{\eco}{\end{corollary}}
\nc{\ere}{\hfill\mbox{$\Diamond$}\end{remark}}
\nc{\eex}{\hfill\mbox{$\Diamond$}\end{example}}
\nc{\epf}{\hfill\mbox{$\Box$}~\\~\\}
\nc{\ot}{\otimes}
\nc{\bsb}{\begin{Sb}}
\nc{\esb}{\end{Sb}}
\nc{\ct}{\mbox{${\cal T}$}}
\nc{\ctb}{\mbox{${\cal T}\sb B$}}
\nc{\bcd}{\beq\begin{CD}}
\nc{\ecd}{\end{CD}\eeq}
\nc{\ba}{\begin{array}}
\nc{\ea}{\end{array}}
\nc{\bea}{\begin{eqnarray}}
\nc{\eea}{\end{eqnarray}}
\nc{\be}{\begin{enumerate}}
\nc{\ee}{\end{enumerate}}
\nc{\beq}{\begin{equation}}
\nc{\eeq}{\end{equation}}
\nc{\bi}{\begin{itemize}}
\nc{\ei}{\end{itemize}}
\nc{\kr}{\mbox{Ker}}
\nc{\te}{\!\ot\!}
\nc{\pf}{\mbox{$P\!\sb F$}}
\nc{\pn}{\mbox{$P\!\sb\nu$}}
\nc{\bmlp}{\mbox{\boldmath$\left(\right.$}}
\nc{\bmrp}{\mbox{\boldmath$\left.\right)$}}
\rnc{\phi}{\mbox{$\varphi$}}
\nc{\LAblp}{\mbox{\LARGE\boldmath$($}}
\nc{\LAbrp}{\mbox{\LARGE\boldmath$)$}}
\nc{\Lblp}{\mbox{\Large\boldmath$($}}
\nc{\Lbrp}{\mbox{\Large\boldmath$)$}}
\nc{\lblp}{\mbox{\large\boldmath$($}}
\nc{\lbrp}{\mbox{\large\boldmath$)$}}
\nc{\blp}{\mbox{\boldmath$($}}
\nc{\brp}{\mbox{\boldmath$)$}}
\nc{\LAlp}{\mbox{\LARGE $($}}
\nc{\LArp}{\mbox{\LARGE $)$}}
\nc{\Llp}{\mbox{\Large $($}}
\nc{\Lrp}{\mbox{\Large $)$}}
\nc{\llp}{\mbox{\large $($}}
\nc{\lrp}{\mbox{\large $)$}}
\nc{\lbc}{\mbox{\Large\boldmath$,$}}
\nc{\lc}{\mbox{\Large$,$}}
\nc{\Lall}{\mbox{\Large$\forall\;$}}
\nc{\bc}{\mbox{\boldmath$,$}}
\rnc{\epsilon}{\varepsilon}
\nc{\ra}{\rightarrow}
\nc{\ci}{\circ}
\nc{\cc}{\!\ci\!}
\nc{\T}{\mbox{\sf T}}
\nc{\can}{\mbox{\em\sf T}\!\sb R}
\nc{\cnl}{$\mbox{\sf T}\!\sb R$}
\nc{\lra}{\longrightarrow}
\nc{\M}{\mbox{\rm Map}}
\rnc{\to}{\mapsto}
\nc{\imp}{\Rightarrow}
\rnc{\iff}{\Leftrightarrow}
\nc{\ob}{\mbox{$\Omega\sp{1}\! (\! B)$}}
\nc{\op}{\mbox{$\Omega\sp{1}\! (\! P)$}}
\nc{\oa}{\mbox{$\Omega\sp{1}\! (\! A)$}}
\nc{\inc}{\mbox{$\,\subseteq\;$}}
\nc{\de}{\mbox{$\Delta$}}
\nc{\spp}{\mbox{${\cal S}{\cal P}(P)$}}
\nc{\dr}{\mbox{$\Delta_{R}$}}
\nc{\dsr}{\mbox{$\Delta_{\Omega^1P}$}}
\nc{\ad}{\mbox{$\mathop{\mbox{\rm Ad}}_R$}}
\nc{\m}{\mbox{m}}
\nc{\0}{\sb{(0)}}
\nc{\1}{\sb{(1)}}
\nc{\2}{\sb{(2)}}
\nc{\3}{\sb{(3)}}
\nc{\4}{\sb{(4)}}
\nc{\5}{\sb{(5)}}
\nc{\6}{\sb{(6)}}
\nc{\7}{\sb{(7)}}
\nc{\hsp}{\hspace*}
\nc{\nin}{\mbox{$n\in\{ 0\}\!\cup\!{\Bbb N}$}}
\nc{\al}{\mbox{$\alpha$}}
\nc{\bet}{\mbox{$\beta$}}
\nc{\ha}{\mbox{$\alpha$}}
\nc{\hb}{\mbox{$\beta$}}
\nc{\hg}{\mbox{$\gamma$}}
\nc{\hd}{\mbox{$\delta$}}
\nc{\he}{\mbox{$\varepsilon$}}
\nc{\hz}{\mbox{$\zeta$}}
\nc{\hs}{\mbox{$\sigma$}}
\nc{\hk}{\mbox{$\kappa$}}
\nc{\hm}{\mbox{$\mu$}}
\nc{\hn}{\mbox{$\nu$}}
\nc{\la}{\mbox{$\lambda$}}
\nc{\hl}{\mbox{$\lambda$}}
\nc{\hG}{\mbox{$\Gamma$}}
\nc{\hD}{\mbox{$\Delta$}}
\nc{\th}{\mbox{$\theta$}}
\nc{\Th}{\mbox{$\Theta$}}
\nc{\ho}{\mbox{$\omega$}}
\nc{\hO}{\mbox{$\Omega$}}
\nc{\hp}{\mbox{$\pi$}}
\nc{\hP}{\mbox{$\Pi$}}
\nc{\bpf}{{\it Proof.~~}}
\nc{\as}{\mbox{$A(S^3\sb s)$}}
\nc{\bs}{\mbox{$A(S^2\sb s)$}}
\nc{\slq}{\mbox{$A(SL\sb q(2))$}}
\nc{\fr}{\mbox{$Fr\llp A(SL(2,\IC))\lrp$}}
\nc{\slc}{\mbox{$A(SL(2,\IC))$}}
\nc{\af}{\mbox{$A(F)$}}
\rnc{\widetilde}{\tilde}
\nc{\suq}{\mbox{$A(SU_q(2))$}}
\nc{\asq}{\mbox{$A(S_q^2)$}}
\nc{\tasq}{\mbox{$\widetilde{A}(S_q^2)$}}

\def\esl{{\mbox{$E\sb{\frak s\frak l (2,{\Bbb C})}$}}}
\def\esu{{\mbox{$E\sb{\frak s\frak u(2)}$}}}
\def\zf{{\mbox{${\Bbb Z}\sb 4$}}}
\def\zt{{\mbox{$2{\Bbb Z}\sb 2$}}}
\def\ox{{\mbox{$\Omega\sp 1\sb{\frak M}X$}}}
\def\oxh{{\mbox{$\Omega\sp 1\sb{\frak M-hor}X$}}}
\def\oxs{{\mbox{$\Omega\sp 1\sb{\frak M-shor}X$}}}
\def\Fr{\mbox{Fr}}
\def\gal{-Galois extension}
\def\hge{Hopf-Galois extension}
\def\ta{\tilde a}
\def\tb{\tilde b}
\def\tc{\tilde c}
\def\td{\tilde d}
\def\st{\stackrel}
\def\d{\mbox{$\mathop{\mbox{\rm d}}$}}
\def\id{\mbox{$\mathop{\mbox{\rm id}}$}}
\def\ker{\mbox{$\mathop{\mbox{\rm Ker$\,$}}$}}
\def\hom{\mbox{$\mathop{\mbox{\rm Hom}}$}}
\def\im{\mbox{$\mathop{\mbox{\rm Im}}$}}
\def\map{\mbox{$\mathop{\mbox{\rm Map}}$}}
\def\o{\sp{[1]}}
\def\t{\sp{[2]}}
\def\mo{\sp{[-1]}}
\def\z{\sp{[0]}}


\newcommand{\Sp}{{\rm Sp}\,}
\newcommand{\Mor}{\mbox{$\rm Mor$}}
\newcommand{\skrA}{{\cal A}}
\newcommand{\Phase}{\mbox{$\rm Phase\,$}}
\newcommand{\diag}{{\rm diag}}
\newcommand{\inv}{{\rm inv}}
\newcommand{\poi}{{\rm pt}}
\newcommand{\Dim}{{\rm dim}\,}
\newcommand{\Ker}{{\rm Ker}\,}
\newcommand{\Mat}{{\rm Mat}\,}
\newcommand{\Rep}{{\rm Rep}\,}
\newcommand{\Fun}{{\rm Fun}\,}
\newcommand{\Tr}{{\rm Tr}\,}
\newcommand{\supp}{\mbox{$\rm supp$}}
\newcommand{\half}{\frac{1}{2}}
\newcommand{\skrF}{{A}}
\newcommand{\skrD}{{\cal D}}
\newcommand{\skrC}{{\cal C}}
\newcommand{\ttimes}{\mbox{$\hspace{.5mm}\bigcirc\hspace{-4.9mm}
\perp\hspace{1mm}$}}
\newcommand{\Ttimes}{\mbox{$\hspace{.5mm}\bigcirc\hspace{-3.7mm}
\raisebox{-.7mm}{$\top$}\hspace{1mm}$}}
\newcommand{\bbr}{{\bf R}}
\newcommand{\bbz}{{\bf Z}}
\newcommand{\Ci}{C_{\infty}}
\newcommand{\Cb}{C_{b}}
\newcommand{\fa}{\forall}
\newcommand{\wrt}{with respect to}
\newcommand{\qg}{quantum group}
\newcommand{\qgs}{quantum groups}
\newcommand{\cs}{classical space}
\newcommand{\qs}{quantum space}
\newcommand{\ch}{character}
\newcommand{\chs}{characters}

\def\inbar{\,\vrule height1.5ex width.4pt depth0pt}
\def\IC{{\Bbb C}}
\def\IR{{\Bbb R}}
\def\IZ{{\Bbb Z}}
\def\IN{{\Bbb N}}
\def\otc{\otimes_{\IC}}
\def\ra{\rightarrow}
\def\ota{\otimes_ A}
\def\otza{\otimes_{ Z(A)}}
\def\otc{\otimes_{\IC}}
\def\h{\rho}
\def\x{\zeta}
\def\th{\theta}
\def\s{\sigma}
\def\<{\langle}
\def\>{\rangle}
\def\vt{\vartriangleright}

\title{
{\Large\bf STRONG CONNECTIONS AND CHERN-CONNES PAIRING IN THE HOPF-GALOIS THEORY
}\vspace*{5mm}\\
}
\author{
{\sc Ludwik D\c{a}browski}
\vspace*{-0mm}\\
\normalsize Scuola Internazionale Superiore di Studi Avanzati
\vspace*{-0mm}\\
\normalsize Via Beirut 2-4, 34014 Trieste, Italy
\vspace*{-0mm}\\
http://www.sissa.it/fm/members/faculty/dabro.html
\vspace*{10mm}\\
{\sc Harald Grosse}
\vspace*{-0mm}\\
\normalsize Institute for Theoretical Physics, University of Vienna
\vspace*{-0mm}\\
\normalsize Boltzmanngasse 5, 1090 Vienna, Austria
\vspace*{-.0mm}\\
http://www.thp.univie.ac.at/local/mathphys
\vspace*{10mm}\\
{\sc Piotr M.~Hajac}
\vspace*{-0mm}\\
\normalsize Mathematical Institute, Polish Academy of Sciences
\vspace*{-0mm}\\
\normalsize  ul.~\'Sniadeckich 8, Warsaw, \mbox{00--950~Poland}
\vspace*{-.0mm}\\
\normalsize and
\vspace*{-.0mm}\\
\normalsize Department of Mathematical Methods in Physics
\vspace*{-0mm}\\
\normalsize Warsaw University, ul.~Ho\.{z}a 74, Warsaw, \mbox{00--682~Poland}
\vspace*{-.0mm}\\
http://www.fuw.edu.pl/$\!\!\!\!\!\!\widetilde{\phantom{mmm}}\!\!\!\!\!\!$pmh
}
\date{}
\maketitle

\begin{abstract}
\noindent\normalsize
We reformulate the concept of connection on a \hge\ $B\inc P$ in order
to apply it in computing the Chern-Connes pairing between the cyclic cohomology
$HC^{2n} (B)$ and $K_0 (B)$.
This reformulation allows us to show that a \hge\ admitting a strong
connection is projective and left faithfully flat.
It also enables us to conclude that a strong connection is a Cuntz-Quillen-type
bimodule connection.
To exemplify the theory, we construct a strong connection
(super Dirac monopole) to find out the Chern-Connes pairing for
the super line bundles associated to super Hopf fibration.
\end{abstract}

\newpage
\section*{Introduction}

A noncommutative-geometric concept of principal bundles and
characteristic classes is given by the Hopf-Galois theory
of ring extensions and the pairing between cyclic cohomology
and $K$-theory, respectively.
In the spirit of the Serre-Swann theorem, the quantum vector
bundles are given as finitely generated projective modules
associated to an $H$\gal\ via a corepresentation of Hopf algebra $H$.
The $K_0$-class of such a module can be then paired with the cohomology
class of a cyclic cocycle to produce an invariant playing the role
of an integrated characteristic class of a vector bundle.
To obtain these invariants,
we provide a theory of connections on \hge s which can be used in calculating
projector matrices of associated quantum vector bundles.
A main point of this paper is that strong connections on
a \hge\ $B\inc P$ are equivalent to left $B$-linear right $H$-colinear
unital splittings of the multiplication map $B\ot P \ra P$.
Since connections can be considered as appropriate liftings
of the translation map (restricted inverse of the canonical
Galois map), knowing a connection yields automatically an
{\em explicit} expression for the translation map.
Vice-versa, an {\em explicit} formula for the translation map
might immediately indicate a formula for connection.
(This is important from the practical point of view.)
If a connection is strong, then  the simple machinery
presented herein helps to extract the projective module data
of an associated quantum vector bundle. One can then
plug it in to the computation of the pairing.
In the classical geometry, characteristic classes of
associated vector bundles are computed from connections on
{\em principal} bundles.
Our approach parallels to some extent this idea in the
quantum-geometric setting.

We work within the general framework of noncommutative geometry,
quantum groups and Galois-type theories.
For an introduction to \hge s we refer to \cite{m-s93, s-hj94} and for
a comprehensive description of the Chern-Connes pairing to
\cite{c-a94,l-jl97}.
The point of view advocated in here was already employed to compute
projector matrices \cite{hm99} and the Chern numbers \cite{h-pm}
of the quantum Hopf line bundles from the Dirac $q$-monopole
connection \cite{bm93}.
Thus, although this work is antedated by \cite{hm99} and \cite{h-pm},
it conceptually precedes these papers, and can be viewed
as a follow up of the theory of connections, strong connections
and  associated quantum vector bundles developed in \cite{bm93},
\cite{h-pm96} and \cite{d-m96}, respectively.
(See \cite[Section~5]{d-m97a} and \cite{d-m97c,d-ma} for an alternative
theory of characteristic classes on  quantum principal bundles.)

We begin in Section~1 by recalling basic facts and definitions.
In Section~2 we first reformulate the concept of general connections
so as to make transparent the characterisation of a strong connection
as an appropriate splitting of the multiplication map $B\ot P \ra P$,
where $P$ is an $H$\gal\ of $B$.
Then we prove the equivalence of four different definitions of
a strong connection, which is the main claim of this paper,
and study its consequences.
As a quick illustration of the theory, we apply it to a strong
and non-strong connection on quantum projective space $\IR P_q^2$.
We obtain, as a by-product, the definition of the ``tangent bundle"
of the Podle\'s equator quantum sphere.
We also show that there are infinitely many canonical strong
connections on the quantum Hopf fibration, and prove that they all
coincide with the Dirac monopole in the classical limit.
A super Dirac monopole is presented in Section~3.
We adapt to the Hopf-Galois setting the construction of a super Hopf
fibration.
Then, employing the super monopole, we compute projector matrices
of the super Hopf line bundles.
Taking advantage of the functoriality of the Chern-Connes pairing,
we conclude that the value of the pairing for the super and classical
Hopf line bundles coincides.
Hence we infer the non-cleftness of the super Hopf fibration.
We end Section~3 by proving that, in analogy with the classical situation,
the direct sum of spin-bundle modules (Dirac spinors) is free of rank two for
both the super and the quantum Hopf fibration.
In Appendix, we complement the four descriptions of
a strong connection by providing (appropriately adapted) four
equivalent actions of gauge transformations on connections.

\section{Preliminaries}

Throughout the paper algebras are assumed to be unital and over a field $k$.
The unadorned tensor product stands for the tensor product over $k$.
Our approach is algebraic, so that we use finite sums.
We use the Sweedler notation $\Delta h=h\1\ot h\2$
(summation understood) and its derivatives.
The letter $S$ and $\he$ signify the antipode and counit, respectively.
The convolution product of two
linear maps from a coalgebra to an algebra is denoted in the following way:
$(f*g)(c):=f(c\1)g(c\2)$.
We use the word ``colinear"
with respect to linear maps that preserve the comodule structure.
(Such maps are also called ``covariant.")
We work with right \hge s and skip writing ``right" for brevity.
For an $H$\gal\ $B\inc P$ we write the canonical Galois isomorphism as
\beq
\chi:=(m\ot\id)\circ (\id\ot\sb B \dr )\, :\; P\ot\sb B P\lra P\ot H \ ,
\eeq
where $\dr : P\ra P\ot H$ stands for the comodule algebra coaction
($\dr p =: p\0 \ot p\1$; again, summation understood),
and $m$ for the multiplication map $P\ot P\ra P$.
We say that a \hge\ is cleft iff there exists a  (unital) convolution-invertible
colinear map $\Phi : H\ra P$, and call $\Phi$ a cleaving map.
The concept of cleftness is close but, as explained in the last paragraph of
\cite[Section~4]{dhs99},
not tantamount to the idea of triviality of a principal bundle.
(Trivial is cleft but not vice-versa.)
A cleaving map is usually assumed to be unital, but since any
non-unital $\Phi$ can be unitalised (e.g., see \cite[p.813]{dt86}
or \cite[Section~1]{hm99}),
this assumption, though technically useful, is conceptually redundant.
It also follows from the defining properties of $\Phi$ that it is injective
(e.g., see \cite[Section~1]{hm99}).

Next, note that the canonical map $\chi$, although cannot be
an algebra homomorphism in general, is always determined by its values on generators.
The same is true for $\chi^{-1}$.
The left $P$-linearity of $\chi^{-1}$ makes it practical to restrict it
from $P\ot H$ to $H$, and define the translation map
\beq\label{tra}
\tau : H\ra P\ot_B P \ , ~~~\tau (h) :=  \chi^{-1} (1\ot h) = : h\o \ot_B h\t
~~~\mbox{(summation understood).}
\eeq
The following are properties of $\tau$ compiled from
\cite{s-hj90b,b-t96}:
 \beq\label{t1}
(\id \ot_B \dr ) \ci \tau = (\tau\ot \id ) \ci \hD \ ,
\eeq
 \beq\label{t1b}
((\mbox{flip}\ci\dr) \ot_B \id )) \ci \tau =
(S\ot\tau)\ci \hD \ ,
\eeq
\beq\label{t2}
\hD_{P\ot_B P}\; \ci \tau = (\tau\ot \id ) \ci \ad \ ,
\eeq
\beq\label{t3}
m \ci \tau = \he \ ,
\eeq
\beq\label{t3b}
\tau (h \tilde h ) = \tilde h\o h\o \ot_B h\t\tilde h\t \ .
\eeq
Here $\hD_{P\ot_B P}$ is the coaction on ${P\ot_B P}$ obtained
via the canonical surjection $\pi_B : P\ot P\ra P\ot_B P$
from the diagonal coaction
\beq\label{dsr}
\hD_{P\ot P} : p\ot p' \longmapsto p\0 \ot {p'}\0 \ot p\1 {p'}\1 \ ,
\eeq
and  $\ad (h):=h\2\ot S(h\1)h\3$ is the right adjoint coaction.

To fix convention and clarify some basic issues,
let us recall that the universal differential calculus $\hO^1A$
 (grade one of the universal differential algebra) can be defined
by the exact sequence
\beq\label{bes}
0\lra\hO^1A\lra A\ot A\lra A\lra 0\, ,
\eeq
i.e., as the kernel of the multiplication map. The
differential is given by $\d a:=1\ot a-a\ot 1$.
We can identify $\hO^1A$ with $A\ot A/k$ as left $A$-modules via the maps
\beq\label{iden}
\hO^1A\ni\sum_i a_i\ot a'_i\to\sum_i a_i\ot \pi_A(a'_i)\in A\ot A/k
\ni x\ot \pi_A(y)\to x\d y\in\hO^1A\, ,
\eeq
where $\pi_A:A\ra A/k$ is the canonical surjection. Similarly, one can identify
$\hO^1A$ with $A/k\ot A$ as right $A$-modules
($\sum_i a_i\ot a'_i\to\sum_i \pi_A(a_i)\ot a'_i$). Consequently, for any left
$A$-module $N$, we have $\hO^1A\ot_AN\cong A/k\ot N$. For any splitting
$\imath: A/k\ra A$ of the canonical surjection ($\pi_A\ci\imath=\id$), we have
an injection $\imath\ot\id:A/k\ot N\ra A\ot N$. Thus there is an injection
\beq
f_\imath:\hO^1A\ot_AN\lra A\ot N\, ,~~~
f_\imath\llp\sum_{i,j} a_{ij}\ot a'_{ij}\ot_A n_j\lrp:=
\sum_{i,j} (\imath\ci\pi_A)(a_{ij})\ot a'_{ij} n_j\, .
\eeq
On the other hand, we have a natural map coming from tensoring (\ref{bes})
on the right with $N$:
\beq
f_N:\hO^1A\ot_AN\lra A\ot N\, ,~~~
f_N\llp\sum_{i,j} a_{ij}\ot a'_{ij}\ot_A n_j\lrp:=
\sum_{i,j} a_{ij}\ot a'_{ij} n_j\, .
\eeq
Since $\pi_A\ci\imath=\id$, we have
$(\pi_A\ot\id)\ci f_N=(\pi_A\ot\id)\ci f_\imath$, whence
\beq
\llp(\imath\ci\pi_A)\ot\id\lrp\ci f_N
=\llp(\imath\ci\pi_A)\ot\id\lrp\ci f_\imath
=f_\imath\, .
\eeq
It follows now from the injectivity of $f_\imath$ that $f_N$ is injective. Thus
we have shown that (\ref{bes}) yields the exact sequence:
\beq
0\lra\hO^1A\ot_AN\lra A\ot N\lra N\lra 0\, .
\eeq
If $B$ is a subalgebra of $P$, then we can also write $(\hO^1B)P$ for the
kernel of the multiplication map $B\ot P\ra P$. Indeed, $m((\hO^1B)P)=0$,
and if $\sum_i b_i\ot p_i\in\ker(B\ot P\st{m}{\ra} P)$, then
\beq
\sum_i b_i\ot p_i
=\sum_i (b_i\ot p_i-1\ot b_i p_i)
=-\sum_i (\d b_i)p_i\in(\hO^1B)P\, .
\eeq
To sum up, we have (cf.~\cite[p.251]{hm99})
\beq\label{id2}
\hO^1B\ot_B P\cong\ker(B\ot P\st{m}{\ra} P)=(\hO^1B)P\, .
\eeq

The following are the universal-differential-calculus
 versions of general-calculus definitions in~\cite{bm93,h-pm96}:
\bde[\cite{bm93}]\label{condef}
Let $B\inc P$ be an $H$\gal. Denote by $\hO^1P$ the universal differential
calculus on~$P$ and by \dsr\ the restriction of $\hD_{P\ot P}$ to  $\hO^1P$.
A left $P$-module projection $\Pi$ on $\hO^1P$ is
called a {\em connection} iff
\beq\label{hf}
\ker\Pi =  P(\Omega\sp{1}B)P ~~\mbox{ (horizontal forms)},
\eeq
\beq\label{rc}
\dsr\ci\Pi = (\Pi\ot \id)\ci\dsr~~\mbox{  (right colinearity)}\, .
\eeq
\ede
\bde[\cite{bm93}] \label{confordef}
Let $P$, $H$, $B$ and $\hO^1P$ be as above.
 A $k$-homomorphism $\omega : H\ra\hO^1P$ such that $\ho(1)=0$
is called a {\em connection form} iff it satisfies the
following properties:
\be
\item $ (m\ot \id)\ci(\id\ot\dr )\ci\omega = 1\ot(\id - \eps)$
(fundamental vector field condition),
\item $\dsr\ci\omega = (\omega\ot \id)\ci \ad$
(right adjoint colinearity).
\ee
\ede
For every \hge\ there is a one-to-one
correspondence between connections and
connection forms (see \cite[p.606]{bm93} or \cite[Proposition~2.1]{m-s97}).
In particular, the connection
$\Pi\sp{\omega}$ associated to
a connection form $\omega$ is given by the
formula:
\beq\label{omefor}
\Pi\sp{\omega}(dp)=p\0\ho(p\1)\, .
\eeq
(Since $\Pi\sp{\omega}$ is a left $P$-module homomorphism,
it suffices to know its values on exact forms.)
\bde[\cite{h-pm96}]\label{scdef}
Let $\Pi$ be a connection in the sense of Definition~\ref{condef}.
It is called {\em strong} iff
$(\id - \Pi)(dP)\inc(\Omega\sp{1}B)P$. We say that a connection form is strong
iff its associated connection is strong.
\ede

Let us now have a closer look at the concept of connection.
For the sake of brevity we put
\beq
\overline{\chi} = (m\ot\id)\ci (\id\ot \dr) : P\ot P\ra P\ot H\, ,
\eeq
 denote by $\widehat{\chi}$ its restriction to $\Omega\sp{1}P$,
and by $H^+$ the kernel of the counit map (augmentation ideal).
Since $\llp (\id\ot\he)\ci \widehat{\chi} \lrp (\Omega\sp{1}P) = 0$,
we have $\widehat{\chi}  (\Omega\sp{1}P)= P\ot H^+$.
Consider $ P\ot H$, and similarly $ P\ot H^+$, as a right comodule via the map
\beq
\hD_{P\ot H} : p\ot h \longmapsto p\0 \ot {h}\2 \ot p\1 {S(h\1)}h\3 \ .
\eeq
Then there is a one-to-one correspondence between connections
and left $P$-linear right $H$-colinear splittings of  $\widehat{\chi}$
\cite[p.606]{bm93}.
Since $H=H^+\oplus k$, we can define
\beq
\overline{\hs}(p\ot h) = \left\{
\begin{array}{ll}
\hs (p\ot h) & \mbox{for $h\in H^+$}\\
p\ot h1_P
& \mbox{for $h \in k$},
\end{array}
\right.
\eeq
where $\hs$ is a splitting of $\widehat{\chi}$.
On the other hand, we can consider unital left $P$-linear right $H$-colinear
splittings $r$ of the canonical surjection $\pi_B : P\ot P \ra P \ot_B P$.
This leads to the following commutative diagram of exact rows
of left $P$-modules right $H$-comodules (see above for the comodule structures):
\bigskip
\beq\label{dia}
\def\normalbaselines{\baselineskip30pt
\lineskip3pt \lineskiplimit3pt }
\def\mapright#1{\smash{
\mathop{\!\!\!-\!\!\!\longrightarrow\!\!\!}
\limits^{#1}}}
\def\mapdown#1{\Big\downarrow
\rlap{$\vcenter{\hbox{$\scriptstyle#1$}}$}}
\def\Mapright#1{\smash{\mathop{-\!\!\!-\!\!\!-\!\!\!-\!\!\!\lra}\limits^{#1}}}
\def\mapleft#1{\smash{\mathop{\longleftarrow\!\!\!-}\limits_{#1}}}
\matrix{
0 &\mapright{} &  P(\hO^1\! B)P & \mapright{} & \hO^1\! P &
\stackrel{\large\Mapright{\widehat{\chi}}}{\mapleft{\sigma}}
&P\ot H^+ &\mapright{} & 0 \cr
&&\bigg\| &&\mapdown{} &&\mapdown{} &&  \cr
0 &\mapright{} &  P(\hO^1\! B)P & \mapright{} & P\ot P &
\stackrel{\large\Mapright{\overline{\chi}}}{\mapleft{\overline{\sigma}}}
&P\ot H &\mapright{} & 0 \cr
&&\bigg\| &&\bigg\| &&\mapdown{\chi^{-1}} &&  \cr
0& \mapright{}   & P(\hO^1\! B)P & \mapright{} & P\ot P &
\stackrel{\large\Mapright{\pi_B}}{\mapleft{r}} & P\ot_B P &\mapright{} &0 \, .\cr
}
\eeq
One can check that
$\overline\chi$ intertwines the relevant comodule structures
\bea
(\hD_{P\ot H}\ci \overline\chi ) (p\ot p')
\!\!\!\!\!\!&&
=\hD_{P\ot H} (p p'\0\ot p'\1 )
\nonumber\\ &&
=p\0 p'\0\ot p'\3 \ot p\1 p'\1 S(p'\2 ) p'\4
\nonumber\\ &&
=p\0 p'\0 \ot p'\1 \ot p\1 p'\2
\nonumber\\ &&
=(\overline\chi\ot \id) (p\0 \ot p'\0\ot p\1 p'\1 )
\nonumber\\ &&
= \llp (\overline\chi\ot \id)\ci \hD_{P\ot P}\lrp (p\ot p') \, .
\eea
Other calculations to verify that this diagram is a commutative
diagram of right $H$-comodules are of the same kind. To see that
$\ker\pi_B = P(\hO^1B)P$ one can argue as above~(\ref{id2}).

Yet another description of a connection as a splitting is as follows.
Denote $\pi_B (\hO^1\! P)$ by $\hO^1_B P$
(relative differential forms as in \cite[Section~2]{cq95}).
The commutativity of the diagram (\ref{dia}) implies that
the restriction of the canonical map $\check \chi : \hO^1_B P\ra P\ot H^+$
is an isomorphism.
Let \ho\ be a connection form and $\widetilde{\ho}$ its restriction to $H^+$.
Similarly, let $\widetilde{\tau}$ be the restriction to $H^+$ of the translation map.
Recall that \hs\ is the left $P$-module extension of $\widetilde{\ho}$
~\cite[p.606]{bm93}.
Hence the commutativity of (\ref{dia})
implies also (for any $\widetilde{\ho}$)
$\pi_B\ci \widetilde{\ho} = \widetilde{\tau}$.
Since the translation map $\tau$ is unital, knowing $\widetilde{\tau}$
is tantamount to knowing $\tau$.
Thus a connection form yields an {\em explicit} expression for the
translation map. On the other hand, viewing $H^+$ as a right comodule
under the right adjoint coaction $\ad$ allows one to define equivalently
a connection as a {\em colinear} lifting of the restricted
translation map $\widetilde{\tau}$.
Indeed, we can complete the equality
$\pi_B\ci \widetilde{\ho} = \widetilde{\tau}$ to the commutative diagram
\vskip.2cm\noindent
\beq\label{dia3}
\def\normalbaselines{\baselineskip30pt
\lineskip3pt \lineskiplimit3pt }
\def\mapright#1{\smash{
\mathop{\!\!\!-\!\!\!\longrightarrow\!\!\!}
\limits^{#1}}}
\def\mapup#1{\Big\uparrow
\rlap{$\vcenter{\hbox{$\scriptstyle#1$}}$}}
\def\mapse#1{\smash{\mathop{\vector(1,-1){20}}\limits^{#1}}}
\matrix{
\hO^1\! P &\mapright{\widehat{\chi}} & P\ot H^+ \cr
\mapup{\widetilde{\omega}}&
\setlength{\unitlength}{1mm}
\begin{picture}(12,0)
\put(2,5){\vector(1,-1){8}}
\put(7,2){\mbox{\scriptsize$ \pi_B$}}
\end{picture}
 &\mapup{\check{\chi}} \cr
H^+ \!\!\!& \mapright{\widetilde{\tau}} & \hO^1_B P
\cr
}
\eeq
\vskip.2cm\noindent
and directly verify this assertion.
This explains the close resemblance between the formulas for the
translation maps and connection forms.
For example, compare (\ref{tau}) with (\ref{omform-}-\ref{omform+}) and
the proof of \cite[Proposition~2.10]{h-pm96} with \cite[2.14]{h-pm96}.
Compare also  \cite[Corollary~2.3]{dhs99} and \cite[Proposition~5.3]{bm93}.

A natural next step is to consider
associated quantum vector bundles. More precisely, what we need here is a
replacement of the module of sections of an associated vector bundle. In the
classical case such sections can be equivalently described  as ``functions
of type $\varrho$" from the total space of a principal bundle to a vector
space. We follow this construction in the quantum case by considering
$B$-bimodules of colinear maps $\mbox{Hom}_\rho(V,P)$ associated with
an $H$\gal\ $B\inc P$ via a corepresentation $\rho:V\ra V\ot H$
(see \cite[Appendix~B]{d-m97a} or \cite{d-m96}).

Proposition~\ref{as} gives a formula for a splitting of the multiplication map
$B\ot \hom_{\rho}(V,P)\ra \hom_{\rho}(V,P)$,
and a splitting of the multiplication map is almost the same as a projector
matrix, for it is an embedding of
$\hom_{\rho}(V,P)$ in the free $B$-module $B\ot \hom_{\rho}(V,P)$.
However, to turn a splitting into
a concrete recipe for producing finite size projector matrices of finitely
generated projective modules, we need the following general lemma:
\ble[\cite{hm99}]\label{genl}
Let $A$ be an algebra and $M$ a projective left $A$-module generated by linearly
independent generators $g_1,...,g_n$.
Also, let $\{\widetilde{g}_\mu\}_{\mu\in I}$
be a completion of $\{g_1,...,g_n\}$ to a linear basis of $M$,
$f_2$ be a left $A$-linear splitting of the multiplication map
$A\ot M\ra M$ given by the formula
$f_2(g_k)=\sum_{l=1}^{n}a_{kl}\ot g_l+\sum_{\mu\in I}a_{k\mu}\ot\widetilde{g}_\mu$,
and  $c_{\mu l}\in A$ a choice of coefficients such that
$\widetilde{g}_\mu=\sum_{l=1}^{n}c_{\mu l}g_l$. Then
$E_{kl}=a_{kl} +\sum_{\mu\in I}a_{k\mu}c_{\mu l}$ defines a projector matrix
of $M$, i.e., $E\in M_n(A)$, $E^2=E$ and $A^n E$ (row times matrix)
and $M$ are isomorphic as left $A$-modules.
\ele

For our later purpose, we also need the following general digression.
Let $A$ be an algebra, and let $E$, $F$ be idempotents in
$M_m(A)$, $M_n(A)$, respectively.
It can be verified that the projective modules $A^m E$ and $A^n F$
are isomorphic if there exist maps $L$ and $\tilde L$
\vspace*{3mm}\beq\label{LL}
\def\normalbaselines{\baselineskip30pt\lineskip3pt\lineskiplimit3pt}
\def\mapright#1{\smash{\mathop{-\!\!\!\lra}\limits^{#1}}}
\def\Mapright#1{\smash{\mathop{-\!\!\!-\!\!\!-\!\!\!-\!\!\!\lra}\limits^{#1}}}
\def\mapleft#1{\smash{\mathop{\longleftarrow\!\!\!-}\limits_{#1}}}
\def\mapdown#1{\Big\downarrow\rlap{$\vcenter{\hbox{$\scriptstyle#1$}}$}}
\def\mapup#1{\Big\uparrow\llap{$\vcenter{\hbox{$\scriptstyle#1$}}$}}
\matrix{
A^m & \stackrel{\large\Mapright{\mbox{\scriptsize $L$}}}{\mapleft{\tilde{L}}} & A^n\cr
 \mapdown{}\mapup{} && \mapdown{}\mapup{} \cr
A^mE & \stackrel{\large\Mapright{}}{\mapleft{}} & A^nF \cr
}
\eeq
such that
\beq\label{ll}
ELF = LF, ~~F\tilde L E = \tilde L E, ~~EL \tilde L = E, ~~F \tilde L L = F .
\eeq

\section{Connections}
\setcounter{equation}{0}

We begin this section by considering general connections on \hge s
as appropriate splittings.
It is known that, under the assumption
of faithful flatness, there always exist a connection on a \hge\
\cite[Satz 6.3.5]{s-p93} (cf.\ \cite[Theorem 4.1]{d-m97a}).
(For a comprehensive review of faithful flatness see \cite{b-n72}.)
Chasing diagram (\ref{dia}) and playing around with appropriate modifications
of its rows we obtain:
\bpr\label{claim0}
Let $B\inc P$ be an $H$\gal.
Denote by ${\cal C}(P)$ the space of connection forms on $P$,
by ${\cal R}(P)$ the space of unital left $P$-linear right $H$-colinear
splittings $r$ of the canonical surjection $\pi_B : P\ot P \ra P \ot_B P$,
and by ${\cal S}(P)$ the space of unital left $B$-linear right $H$-colinear
maps $s : P\ra P\ot P  $ satisfying
$(\pi_B \ci s )(p) = 1\ot_B p$.
Then the formulas
\beq
\Psi ( \ho ) \, (p\ot_B p') = pp'\ot 1 + pp'\0 \ho(p'\1 ), ~~~
\widetilde{\Psi}(r) \, (h) = (r\ci \tau ) (h -\he (h)),
\eeq
define mutually inverse bijections
$
{\cal C}(P)
\stackrel{\Psi}{\ra}
{\cal R}(P)
\stackrel{\widetilde\Psi}{\ra}
{\cal C}(P)
$
and, similarly, the formulas
\beq
\Xi ( r ) \, (p) = r(1\ot_B p), ~~~
\widetilde{\Xi}(s) \, (p\ot_B p') = ps(p').
\eeq
determine mutually inverse bijections
$
{\cal R}(P)
\stackrel{\Xi}{\ra}
{\cal S}(P)
\stackrel{\widetilde\Xi}{\ra}
{\cal R}(P).
$
\epr
\bpf
Let us first check that $\Psi({\cal C}(P))\inc {\cal R}(P)$.
It is clear that, for any $\ho\in {\cal C}(P)$,
$\Psi(\ho)$ is unital and left $P$-linear.
To see that $\Psi(\ho)$ is $H$-colinear, we use (\ref{omefor}) and (\ref{rc}):
\bea
\llp \hD_{P\ot P}\ci \Psi(\ho)\lrp (p\ot_B p')
\!\!\!\!\!\!&&
= \hD_{P\ot P}(p p'\ot 1 + \Pi^\omega (p \d p'))
\nonumber\\ &&
= (p\0 p'\0\ot 1 + \Pi^\omega (p\0 \d p'\0))\ot p\1 p'\1
\nonumber\\ &&
= \Psi(\ho)(p\0 \ot_B p'\0)\ot p\1 p'\1
\nonumber\\ &&
= \llp(\Psi(\ho)\ot \id) \ci \hD_{P\ot_B P}\lrp (p\ot_B p') \, .
\eea
To verify that  $\Psi(\ho)$ is a splitting of the canonical
surjection $\pi_B$, recall that $\ker \pi_B =  P(\hO^1\! B)P$ and note that
$(\Pi^\omega)^2 = \Pi^\omega$ entails
$\ker \Pi^\omega =  (\id - \Pi^\omega ) (\hO^1\! P)$. Thus, by (\ref{hf}),
we have $\ker \pi_B = (\id - \Pi^\omega ) (\hO^1\! P)$. Combining this with
(\ref{omefor}) we obtain
\bea
\llp\id - \pi_B\ci \Psi(\ho) \lrp  (p\ot_B p')
\!\!\!\!\!\!&&
= \pi_B \llp p\ot p' - pp'\ot 1 - \Pi^\omega (p\d p')\lrp
\nonumber\\ &&
= \llp \pi_B \ci (\id - \Pi^\omega ) \lrp (p\d p') = 0 \, .
\eea

The next step is to check that
 $\widetilde\Psi({\cal R}(P))\inc {\cal C}(P)$.
To see that $\widetilde\Psi(r)\, (H) \inc \hO^1\! P$ for any $r\in {\cal R}(P)$,
we take advantage of  property (\ref{t3}) of the translation map $\tau$,
and compute:
\beq
\llp m \ci \widetilde\Psi(r) \lrp h =
\llp m \ci \pi_B \ci r \ci \tau \lrp (h-\he(h)) =
\llp m \ci \tau \lrp (h-\he(h)) =
\he (h-\he(h)) = 0 \, .
\eeq
(Here we abuse the notation and denote also by $m$ the multiplication map on
$P\ot_B P$.)
It is immediate that $\widetilde\Psi(r)(1) = 0$.
Furthermore, using the colinearity of $r$ and  (\ref{t3}),
we verify the colinearity of $\widetilde\Psi(r)$:
$ \dsr\ci\widetilde\Psi(r) (h) = (\ho\ot\id)\ci \ad $.
To check the fundamental-vector-field condition we note
\beq
 (\overline{\chi}\ci r\ci \tau ) (h-\he(h)) = 1\ot (h-\he(h))\, ,
\eeq
which is equivalent to $
 (\overline{\chi}\ci r\ci \chi^{-1}) = \id $.
The latter equality, however, follows from the commutativity of
(\ref{dia}), as needed.

It remains to show that
$\widetilde\Psi\ci\Psi = \id$ and
$\Psi\ci\widetilde\Psi = \id$.
To this end, taking advantage of the unitality of $\Psi(\ho)$,
(\ref{t1}) and (\ref{t3}), we compute:
\bea
\llp(\widetilde\Psi\ci\Psi )(\ho) \lrp  (h)
\!\!\!\!\!\!&&
= (\Psi(\ho)\ci\tau) (h- \he (h))
\nonumber\\ &&
= \Psi(\ho)(h\o\ot_B h\t) - \he (h)\ot 1
\nonumber\\ &&
= \he (h)\ot 1 + h\o {h\t}\0 \ho ({h\t}\1 )- \he (h)\ot 1
\nonumber\\ &&
= {h\1}\o {h\1}\t \ho ({h\2})
\nonumber\\ &&
= \he({h\1}) \ho ({h\2}) =  \ho ({h})\ .
\eea
Similarly, taking advantage of the unitality and left $P$-linearity of $r$,
we compute
\bea
\llp(\Psi\ci\widetilde\Psi )(r) \lrp  (p\ot_B p')
\!\!\!\!\!\!&&
=pp'\ot 1 + pp'\0 \llp\widetilde\Psi (r)\lrp (p'\1)
\nonumber\\ &&
=pp'\ot 1 + pp'\0 (r\ci \tau) (p'\1 - \he({p'\1}) )
\nonumber\\ &&
=pp'\ot 1 + r\llp pp'\0\tau (p'\1)\lrp - pp'\0 \he({p'\1})\ot 1
\nonumber\\ &&
= r \llp \chi^{-1}( \chi(p\ot_B p'))\lrp = r (p\ot_B p') \ .
\eea
Finally, the proof concerning $\Xi$ and $\widetilde\Xi$ is straightforward.
\epf
\vspace{-.5cm}
\bco\label{00}
An $H$\gal\ $B\inc P$ admits a connection, if there exists
a (not necessarily unital) left $B$-linear right $H$-colinear
map $s : P\ra P\ot P  $ satisfying
$(\pi_B \ci s )(p) = 1\ot_B p$.
\eco
\bpf
Denote by $\overline{{\cal S}(P)}$ the space of all maps $s$
defined in the corollary.
To construct a ``unitalising'' map
$ {\cal T} : \overline{{\cal S}(P)} \ra {\cal S}(P) $,
note that, for $\overline s \in \overline{{\cal S}(P)}$
and $s\in {\cal S}(P) $,
we must have
$s(1) -\overline s (1) = 1\ot 1 -\overline s (1)$.
Extending this equality by the left $P$-linearity,
we can define
\beq
{\cal T}\! (\overline s) (p) = \overline s(p) + p(1\ot 1 - \overline s(1)).
\eeq
It is straightforward to check that
${\cal T}(\overline{{\cal S}(P)})\inc {\cal S}(P)$,
as needed.
\epf
When we think of a connection as an element $s\in {\cal S}(P)$,
then the strongness condition (see Definition~\ref{scdef})
can be put as $s(P) \inc B\ot P$.
(Shift the second term on the right hand side to the left hand side in
\cite[(11)]{m-s97}.)
Describing strong connections as strong elements in ${\cal S}(P)$
is the main point of the below theorem.
The second description is in terms of a covariant differential,
and was hinted at in \cite[Remark~4.3]{h-pm96}.
The third one coincides with the definition of a strong connection
except that we change the inclusion condition
$(\id - \Pi ) (\d P) \inc (\hO^1\! B) P$
to the equivalent equality condition
$(\id - \Pi ) (B\d P) = (\hO^1\! B) P$.
The last description is precisely the definition of a strong connection form.
Thus there is no essentially new approach to connections in the following
theorem. However, proving the equivalence of a strong connection to an
appropriate splitting of the multiplication map $B\ot P \ra P$
enables us to derive several desirable consequences.
We write everything explicitly so as to provide a self-contained and coherent
treatment of the strong connection.
\bth\label{claim1}
Let $B\inc P$ be an $H$\gal. The following are equivalent descriptions
of a strong connection:
\vspace*{2.5mm}\\
1) A unital left $B$-linear right $H$-colinear splitting $s$ of the multiplication
map
$
\def\mapl#1{\smash{\mathop{\leftarrow}\limits_{#1}}}
\def\Mapr#1{\smash{\mathop{-\!\!\!\ra}\limits^{#1}}}
B\ot P
\st{\large\Mapr{m}}{\mapl{s}}
P
$.
\vspace*{2.5mm}\\
2) A right $H$-colinear homomorphism $D: P\ra(\hO^1\! B)P$
annihilating 1 and satisfying the Leibniz rule:
$D(bp)= bDp + \d b.p,\,\fa\, b\in B,\, p\in P$.
\vspace*{2.5mm}\\
3) A left $P$-linear right $H$-colinear projection $\Pi:\hO^1\! P\ra\hO^1\! P$
($\Pi^2=\Pi$) such that\linebreak
 $(\id -\Pi)(B\d P)=(\hO^1\! B)P$.
\vspace*{2.5mm}\\
4) A homomorphism $\ho: H\ra\hO^1\! P$ vanishing on 1 and satisfying:
\vspace*{-2.5mm}
\bi
\item[a)] $
 \dsr\ci\ho=(\ho\ot\id)\ci \ad
$\item[b)] $
 (m\ot\id)\ci(\id\ot\dr)\ci\ho=1\ot(\id-\he)
$\item[c)] $
  \d p-p\0\ho(p\1)\in (\hO^1\! B)P,\;\fa\, p\in P .
$\ei
\ethe
\bpf
Let $V_i\, ,i\in\{1,2,3,4\}$, denote the corresponding spaces of homomorphisms
defined in points 1)--4). We need to construct 4 mappings
\beq
J_1:V_1\ra V_2,\;\;\; J_2:V_2\ra V_3,\;\;\; J_3:V_3\ra V_4,\;\;\; J_4:V_4\ra V_1,
\eeq
satisfying 4 identities:
\beq\label{cycl}
J_4\ci J_3\ci J_2\ci J_1 =\id \mbox{~ and cyclicly permuted versions.}
\eeq

Put $J_1(s)(p)=1\ot p-s(p)$.
(Compare with the right-handed version \cite[(55)]{cq95}.)
Evidently, $J_1(s)$ is a right $H$-colinear
homomorphism from $P$ to $(\hO^1\! B)P=\ker(m:B\ot P\ra P)$  (see (\ref{id2}))
annihilating~1. As for the Leibniz rule, we have
\beq
J_1(s)(bp)=1\ot bp-s(bp)=\d b.p+b\ot p-bs(p)=\d b.p +bJ_1(s)(p).
\eeq
This establishes $J_1$ as a map from $V_1$ to $V_2\,$.

Next, put $J_2(D)(p'\d p)=p'(\d -D)(p)$. Observe first that $J_2(D)$
is a well-defined endomorphism of $\hO^1\! P$ because $D1=0$ (see (\ref{iden})).
Choose $b_i\in B$, $p_i\in P$, such that 
$Dp=\sum_i (\d b_i)p_i=\sum_i (\d (b_ip_i)-b_i \d p_i)$. It follows from the
Leibniz rule that
 $J_2(D)\ci\d$ is a left $B$-module map. Thus we have:
\bea
J_2(D)^2(p'\d p)
\!\!\!\!\!\!&&
=J_2(D)(p'\d p)-J_2(D)(p'Dp)
\nonumber\\ &&
=J_2(D)(p'\d p)-\mbox{$\sum_i$}\, p' J_2(D)(\d b_i.p_i)
\nonumber\\ &&
=J_2(D)(p'\d p)-\mbox{$\sum_i$}\, p'\llp J_2(D)\ci\d\lrp(b_i p_i)
+\mbox{$\sum_i$}\, p' b_i\llp J_2(D)\ci\d\lrp(p_i)
\nonumber\\ &&
=J_2(D)(p'\d p),
\eea
Hence $J_2(D)$ is a projection.
Furthermore, note that for any $b_i\in B,\, p_i\in P$, we have
\beq
(\id- J_2(D))(\mbox{$\sum_i$}b_i.\d p_i)
=\mbox{$\sum_i$}b_i.Dp_i\in(\hO^1\! B)P,
\eeq
i.e., $(\id- J_2(D))(\mbox{$\sum_i$}B\d P)\inc(\hO^1\! B)P$.
To see the reverse inclusion, take any
$\mbox{$\sum_i$}b_i\ot p_i\in\ker(B\ot P\st{m}{\ra}P) =(\hO^1\! B)P$.
Then, using the above calculation and the Leibniz rule, we obtain
\beq
0
=\mbox{$\sum_i$}D(b_ip_i)
=\mbox{$\sum_i$}b_iDp_i+\mbox{$\sum_i$}\d b_i.p_i
=(\id -J_2(D))(\mbox{$\sum_i$}b_i.\d p_i)-\mbox{$\sum_i$}b_i\ot p_i\, ,
\eeq
i.e., $\mbox{$\sum_i$}b_i\ot p_i\in\im(\id-J_2(D))$, as needed.

To construct $J_3$, note first that $(\Pi\ci\d) :P\ra\hO^1\! P$ is left $B$-linear.
Indeed, since $\Pi^2=\Pi$, the condition $(\id-\Pi)(B\d P)=(\hO^1\! B)P$
entails $\Pi((\hO^1\! B)P)=0$. Consequently
\beq
\Pi\d(bp)=\Pi(\d b.p)+\Pi(b\d p)=b(\Pi\ci\d)(p),
\eeq
as claimed. Therefore it makes sense to put $ J_3(\Pi)(h)=h\o\Pi(\d h\t)$
(see (\ref{tra})).
This formula defines a homomorphism from $H$ to $\hO^1\! P$ vanishing on~1.
Furthermore, by the right $H$-colinearity of $\Pi$ and property (\ref{t2})
of the translation map, we have
\bea
\llp\dsr\ci J_3(\Pi)\lrp(h)
\!\!\!\!\!\!&&
={h\o}\0\Pi(\d {h\t}\0)\ot {h\o}\1 {h\t}\1
\nonumber\\ &&
={h\2}\o\Pi(\d {h\2}\t)\ot S(h\1)h\3
\nonumber\\ &&
=\llp(J_3(\Pi)\ot\id)\ci \ad\lrp(h).
\eea
As for the property b), note first that $(\id - \Pi ) (B\d P) = (\hO^1\! B)P$
implies, by the left $P$-linearity of $\Pi$, that
$(\id - \Pi )(\hO^1\! P )= P(\hO^1\! B)P$.
Secondly, recall that
$ \overline{\chi} ( P(\hO^1\! B)P) = 0$ (see (\ref{dia})).
Hence
\bea
\llp  \overline{\chi} \ci J_3(\Pi) \lrp (h)
\!\!\!\!\!\!&&
={h\o} \llp(m\ot \id ) \ci (\id \ot \dr )\lrp
\llp \Pi\d {h\t} + (\id - \Pi ) \d {h\t}\lrp
\nonumber\\ &&
={h\o} \llp(m\ot \id ) \ci (\id \ot \dr )\lrp ( 1\ot h\t -h\t \ot 1 )
\nonumber\\ &&
={h\o} {h\t}\0 \ot {h\t}\1 - {h\o} {h\t} \ot 1
\nonumber\\ &&
=1\ot (h -\eps(h)).
\eea
To verify c), we compute:
\beq
\d p - p\0 {p\1}\o \Pi \d {p\1}\t =
\d p - \Pi \d p = (\id - \Pi )\d p \in (\hO^1\! B)P.
\eeq
Consequently, $J_3$ is a mapping from $V_3$ to $V_4$.

Finally, put
\beq\label{sf}
J_4(\ho )(p)= p\ot 1 + p\0 \ho (p\1 ) \ .
\eeq
To see that $ J_4(\ho )$ takes values in $B\ot P$, note that
\bea
p\ot 1 + p\0 \ho (p\1 )
\!\!\!\!\!\!&&
=p\ot 1 -1\ot p +1\ot p + p\0 \ho (p\1 )
\nonumber\\ &&
= 1\ot p - \llp \d p - p\0 \ho (p\1 ) \lrp \in B\ot P
\eea
by property c) of $\ho$. The right $H$-colinearity of $ J_4(\ho )$
follows from property a) of $\ho$.
The remaining needed properties of $ J_4(\ho )$ are immediate.
Consequently, $ J_4$ is a mapping from $V_4$ to $V_1$.

To end the proof, we need to show
$ J_4 \ci J_3 \ci J_2 \ci J_1 = \id$ and its three
cyclicly permuted versions.
We use recurrently the fact that the translation map $\tau$ provides
the inverse of the canonical map $\chi$, so that
$p\0 {p\1}\o \ot_B {p\1}\t = 1\ot_B p$ and $h\o {h\t}\0 \ot {h\t}\1 = 1\ot h$.
\bea
\llp J_4 \ci J_3 \ci J_2 \ci J_1\lrp \! (s) (p)
\!\!\!\!\!\!&&
=p\ot 1 + p\0 \llp  J_3 \ci J_2 \ci J_1 \lrp\! (s) (p\1 )
\nonumber\\ &&
=p\ot 1 + p\0 {p\1}\o \llp J_2 \ci J_1 \lrp\! (s) \llp  \d {p\1}\t \lrp
\nonumber\\ &&
=p\ot 1 + \llp J_2 \ci J_1 \lrp\! (s) (  \d p)
\nonumber\\ &&
=p\ot 1 + \d p - J_1 (s) ( p)
\nonumber\\ &&
=1\ot p - 1\ot p + s(p) = s(p) \ ,
\eea
\bea
\llp J_3 \ci J_2 \ci J_1 \ci  J_4 \lrp \! (\ho ) (h)
\!\!\!\!\!\!&&
=h\o \llp J_2 \ci J_1 \ci  J_4 \lrp \! (\ho ) (\d h\t )
\nonumber\\ &&
=h\o \llp \d - \llp J_1 \ci  J_4 \lrp \! (\ho ) \lrp (h\t )
\nonumber\\ &&
=h\o \llp \d - 1\ot \id + J_4 \!(\ho ) \lrp (h\t )
\nonumber\\ &&
=h\o \llp J_4 \! (\ho ) - \id\ot 1\lrp (h\t )
\nonumber\\ &&
= h\o {h\t}\0 \ho ({h\t}\1) = \ho (h) \ ,
\eea
\bea
\llp J_2 \ci J_1 \ci J_4 \ci J_3\lrp \! (\Pi) (\d p)
\!\!\!\!\!\!&&
=\d p - \llp J_1 \ci J_4 \ci J_3 \lrp\! (\Pi) ( p)
\nonumber\\ &&
=1\ot p - p\ot 1 - 1\ot p + \llp J_4 \ci J_3 \lrp (\Pi) ( p)
\nonumber\\ &&
= -p\ot 1 + p\ot 1  + p\0 \llp J_3 (\Pi) \lrp (p\1 )
\nonumber\\ &&
=p\0 {p\1}\o \Pi (\d {p\1}\t) = \Pi (\d p) \ ,
\eea
\bea
\llp J_1 \ci J_4 \ci J_3 \ci J_2\lrp \! (D) (p)
\!\!\!\!\!\!&&
=1\ot p  - \llp J_4 \ci J_3 \ci J_2 \lrp\! (D) (p)
\nonumber\\ &&
=\d p - p\0 \llp J_3 \ci J_2 \lrp (D) (p\1 )
\nonumber\\ &&
=\d p - p\0 {p\1}\o J_2 (D) (\d {p\1}\t )
\nonumber\\ &&
=\d p - J_2 (D) (\d p) = Dp \ .
\eea
This shows that the maps $J_i$ are bijective.
\epf
\vspace{-1cm}
\bco\label{ff}
If $B\inc P$ is an $H$\gal\ admitting a strong connection, then
\vspace*{2.5mm}\\
1) $P $ is projective as a left $B$-module,
\vspace*{2.5mm}\\
2) $B$ is a direct summand of $P$ as a left $B$-module,
\vspace*{2.5mm}\\
3) $P$ is left faithfully flat over $B$.
\eco
\bpf
Let $s : P\ra B\ot P$ be the splitting associated to a strong connection.
Due to the unitality of $B$ the multiplication map $B\ot P\ra P$ is surjective.
Thus $P$ is  a direct  summand of $B\ot P$
via $s$, and the projectivity of $P$ follows from the freeness of $B\ot P$.

Let $f_P$ be a unital linear functional on $P$.
Then $(\id\ot f_P)\ci s $ is a left $B$ linear map
splitting the inclusion $B\inc P$.
Hence $B$ is a direct summand of $P$.

Finally, since $P$ is projective it is flat.
On the other hand, since $P$ contains $B$ as a direct summand,
it is also faithfully flat.
\epf
In fact, since $s$ embeds $P$ in $B\ot P$ colinearly,
we can say that $P$ is an $H$-equivariantly projective  left $B$-module.
Next, we translate $s$ to the setting of associated quantum bundles so as
to be able to compute their projector matrices with the help of Lemma~\ref{genl}.
\bpr\label{as}
Let $s: P\ra B\ot P$ be the splitting associated to a strong connection on $H$\gal\
$B\inc P$, and let $\rho : V \ra V\ot H$ be any finite dimensional
corepresentation of $H$.
Denote by $\overline{\ell}$ the canonical isomorphism
$B\ot \hom (V, P) \ra \hom(V, B\ot P)$.
Then the formula
\beq\label{sro}
s_{\rho} (\xi ) = \overline{\ell}^{-1} (s\ci \xi)
\eeq
gives a left $B$-linear splitting of
the multiplication map $B\ot \hom_{\rho} (V, P) \ra \hom_{\rho} (V, P) $.
\epr
\bpf
Note first that, since $s$ is right colinear,
$s_{\rho} \llp\hom_{\rho} (V, P) \lrp\inc \hom_{\rho} (V, B\ot P) $.
We need to show that
$\overline{\ell} \llp B\ot \hom_{\rho} (V, P) \lrp
= \hom_{\rho} (V, B\ot P) $,
where $\overline{\ell} ( b\ot \xi ) (v) = b\ot \xi (v)$.
For this purpose we can reason as in the proof of \cite[Proposition~2.3]{hm99}
and construct the following
commutative diagram with exact rows:
\beq\label{5iso}
\def\normalbaselines{\baselineskip30pt
\lineskip3pt \lineskiplimit3pt }
\def\mapright#1{\smash{
\mathop{\!\!\!-\!\!\!\longrightarrow\!\!\!}
\limits^{#1}}}
\def\mapdown#1{\Big\downarrow
\rlap{$\vcenter{\hbox{$\scriptstyle#1$}}$}}
\matrix{
\! 0\!&\mapright{}&\! B\ot\hom_\rho (V, P)\!
&\mapright{}&\! B\ot\hom (V, P)\!
&\mapright{id\otimes \underline{\rho}}&
\! B\ot\hom (V, P\ot H)\!
\!\!\!\!\!\!\!\!\!\!\!\!\!\!\cr
&&\mapdown{}&&\mapdown{\overline{\ell}}&&\mapdown{\underline{\ell}}\cr
\! 0\!&\mapright{}&\! \hom_\rho (V, B\ot P)\!
&\mapright{}&\! \hom (V, B\ot P)\!
&\mapright{\overline{\rho}}&\! \hom (V, B\ot P\ot H)\ .
\!\!\!\!\!\!\!\!\!\!\!\!\!\!\cr
}
\eeq\ \\
Here $\underline{\rho}$ is defined by
$\underline{\rho} (\xi ) = (\xi\ot \id )\ci \rho - \dr \ci \xi$,
and similarly $\overline{\rho}$.
The map $\underline{\ell}$ is the appropriate canonical isomorphism.
Completing the diagram to the left with zeroes and applying the
Five Isomorphism Lemma shows that the restriction  of $\overline{\ell}$
to $B\ot\hom_\rho (V, P)$ is an isomorphism onto $\hom_\rho (V, B\ot P)$,
as needed. Thus $s_\rho$ is a map from  $\hom_\rho (V, P)$ to $B\ot\hom_\rho (V, P)$,
as claimed. Explicitly, $\overline{\ell}^{-1}$ is given by
\beq
\overline{\ell}^{-1} (\varphi )= \sum_i \varphi (e_i ) e^i
= \sum_i \varphi (e_i )\mo \ot\varphi (e_i )\z e^i ,
\eeq
where $\{ e_i \}$ is a basis of $V$, $\{ e^i \}$ its dual, and
we put
$\varphi(v )= \varphi (v )\mo \ot\varphi (v )\z $
(summation understood). Similarly, we can write
$s_\rho (\xi) = s_\rho (\xi )\mo \ot s_\rho (\xi )\z $.)
The left $B$-linearity of $s_\rho$ follows from the
left $B$-linearity of $s$ and $\overline{\ell}$.
Finally, $s_\rho$ splits the multiplication map
because $s$ splits the multiplication map:
\bea
(m\ci s_\rho ) (\xi)(v)
\!\!\!\!\!\!&&
= s_\rho (\xi)\mo s_\rho (\xi)\z (v)
\nonumber\\ &&
= m\Llp \overline{\ell} \llp s_\rho (\xi) \lrp (v) \Lrp
\nonumber\\ &&
= m \llp (s\ci \xi )(v) \lrp
\nonumber\\ &&
= (m \ci s\ci \xi )(v)
\nonumber\\ &&
= \xi (v) \ .
\eea
\epf
Applying the standard reasoning as used in the proof of Corollary~\ref{ff}, we can infer
(under the assumptions of Proposition~\ref{as}) that $\hom_\rho (V, P)$ is
projective as a left $B$-module. On the other hand, if $P$ is left faithfully
flat over $B$ and the antipode of $H$ is bijective, one can prove that
$\hom_\rho (V, P)$ is finitely generated as a left $B$-module~\cite{s-p}.
Thus point 3) of Corollary~\ref{ff} leads to the following conclusion
(cf.\ \cite[Appendix B]{d-m97a}):
\bco\label{proco}
Let $H$ be a Hopf algebra with a bijective antipode,
$B\inc P$ an $H$\gal\ admitting a strong connection,
and $\rho : V\ra V\ot H$ a finite-dimensional corepresentation of $H$.
Then the associated module of colinear maps
$\hom_\rho (V, P)$ is finitely generated projective as a left $B$-module.
\eco
Closely related to $B$-bimodule $\hom_\rho (V, P)$ is $B$-bimodule $P_{\rho} :=
\sum_{\varphi\in \mbox{\scriptsize\rm Hom}_\rho (V, P)} \varphi (V)\inc P$
(cf.\ \cite[Appendix B]{d-m97a}). It turns out that such submodules of $P$
are invariant under the splitting associated to a strong connection:
\bpr\label{pero}
Let $s$ be the splitting associated to a strong connection on an $H$\gal\
$B\inc P$.
Let $\rho : V\ra V\ot H$ be a finite-dimensional corepresentation of $H$
and $P_{\rho} :=
\sum_{\varphi\in \mbox{\scriptsize\rm Hom}_\rho (V, P)} \varphi (V)$.
Then $s (P_{\rho}) \inc B\ot P_{\rho}$.
\epr
\bpf
If $p\in P_{\rho}$ then there exists finitely many
$\widetilde{\varphi}_\nu\in \hom_\rho (V, P)$ such that
\beq\label{pero0}
p = \sum_{\nu} \widetilde{\varphi}_\nu (v_\nu )=
\sum_{\nu} \sum_{k=1}^{\mbox{\scriptsize\rm dim}V} v_{\nu k}
\widetilde{\varphi}_\nu (e_k)
=\sum_{k=1}^{\mbox{\scriptsize\rm dim}V} {\varphi}_k (e_k).
\eeq
Here $\{e_k \}$ is a basis of $V$ and
${\varphi}_k := \sum_{\nu} v_{\nu k}
\widetilde{\varphi}_\nu $.
(Since  $v_{\nu k}$ are simply the coefficients of $v_{\nu}$
with respect to $\{e_k \}$, we have ${\varphi}_k \in \hom_\rho (V, P)$.)
Next, we can always write $s(p) = \sum_{\mu} f_{\mu}\ot s(p)_\mu$,
where $\{ f_\mu \}$ is a linear basis of $B$.
(We have the strongness condition $s (P)\inc B\ot P $.)
Since $s$ and ${\varphi}_k$ are both colinear, so is their composition
$s\ci {\varphi}_k$, and we have
\beq\label{pero1}
\hD_{P\ot P} (s\ci {\varphi}_k )(e_\ell )=
\sum_{m=1}^{\mbox{\scriptsize\rm dim}V} s({\varphi}_k (e_m)) \ot u^\rho_{m \ell}=
\sum_{m=1}^{\mbox{\scriptsize\rm dim}V}
\sum_{\mu} f_\mu \ot  (\varphi_k (e_m))_\mu \ot u^\rho_{m \ell}\ ,
\eeq
where $u^\rho_{m \ell}$ are the matrix elements of corepresentation $\rho$.
On the other hand, remembering that $s (P)\inc B\ot P $, we have
\beq\label{pero2}
\hD_{P\ot P} (s\ci {\varphi}_k )(e_\ell )=
\sum_{\mu} f_\mu \ot \dr\llp{\varphi}_k (e_\ell)\lrp_\mu\ .
\eeq
Combining the above two equalities and
 using the linear independence of $f_\mu$, we obtain
\beq
\dr\llp{\varphi}_k (e_\ell)\lrp_\mu =
\sum_{m=1}^{\mbox{\scriptsize\rm dim}V}
(\varphi_k (e_m))_\mu \ot u^\rho_{m \ell}\ .
\eeq
Hence we can define a bi-index family of $\rho$-colinear maps by the equality
$\varphi_{k \mu}  (e_\ell ) = (\varphi_k (e_\ell ))_\mu$.
Consequently, due to (\ref{pero0}), we have
\beq
s(p) =
s \llp \sum_{k=1}^{\mbox{\scriptsize\rm dim}V} {\varphi}_k (e_k)\lrp
=\sum_{k=1}^{\mbox{\scriptsize\rm dim}V}
\sum_{\mu} f_\mu \ot  (\varphi_k (e_k))_\mu
= \sum_{\mu} f_\mu \ot \llp \sum_{k=1}^{\mbox{\scriptsize\rm dim}V}
\varphi_{k \mu} (e_k)\lrp \in B\ot P_\rho \, ,
\eeq
as claimed.
\epf
\vspace*{-1cm}
\bre\em
Just as we define $J_1$ in the proof of Theorem~\ref{claim1},
we can define the covariant derivative on  $\hom_\rho (V, P)$
via the formula
\beq\label{cd}
\nabla : \hom_\rho (V, P)\ra \hO^1\! B \ot_B \hom_\rho (V, P)~, ~~~
\nabla \xi = 1\ot \xi - s_\rho (\xi )\ .
\eeq
Using identifications in Theorem~\ref{claim1} (isomorphisms $J_i$),
one can check that (\ref{cd}) agrees with \cite[(2.2)]{hm99}.
\ere
\bre\label{cq}
\em
Assume that $\overline s\in\overline{{\cal S}(P)}$ (see Corollary~\ref{00})
enjoys in addition the strongness property
$\overline s (P)\inc B\ot P$.
Then the space of connections ${\cal T}(\overline s)$
can be fully characterized
as the set of all connections $s$ satisfying
\beq
\exists\; \hb\in \hO^1 B~~\forall\; p\in P ~:~~ s(p)-p\hb\in B\ot P .
\eeq
Indeed, one can check that $1\ot 1 - \overline s (1) \in \hO^1 B$,
and the rest follows from the formula for ${\cal T }$.
It is tempting to call such connections semi-strong.
\ere

We now proceed to establishing a link between strong connections and
Cuntz-Quillen connections on bimodules \cite[p.283]{cq95}.
Let $C$ be a coalgebra and $N_1, \, N_2$ right $C$-comodules.
Denote by $A := \hom (C, k) $ the algebra dual to~$C$.
Then $N_1$ and $N_2$ enjoy the following natural
left $A$-module structure (e.g., see \cite[Section~1.6]{m-s93}):
\beq
A\ot N_i \ni a \ot n \to n\0 a(n\1 ) \in N_i \, ,~~~i\in\{ 1,2\} .
\eeq
With respect to this structure, any $k$-homomorphism
from $N_1$ to $ N_2$ is right $C$-colinear if and only
if it is right $A^{op}$-linear. Thus, for an $H$\gal\ $B\inc P$,
algebra $P$ is a $(B, (H^{*})^{op})$-bimodule,
where $H^* := {\rm Hom} (H, k) $ is the algebra dual to $H$ considered
as a coalgebra.
By Theorem~\ref{claim1} (point 2), a strong connection can be given as a
right $(H^*)^{op}$-linear map $D:P\ra \hO^1B\ot_B P$ (see (\ref{id2}))
satisfying the left Leibniz rule and vanishing on~1. Therefore
it seems natural to generalize
the concept of a left bimodule connection \cite[p.284]{cq95} to
\bde\label{lbicon}
Let $N$ be an $(A_1, A_2)$-bimodule. We say that
$
\nabla_L : N \ra \hO^1 A_1 \ot_{A_1} N
$
is a {\em left bimodule connection} iff
it is right $A_2$-linear and satisfies the left Leibniz rule:
$\nabla_L(an)=a\nabla_L(n)+\d a\ot_{A_1} n$, $\;\forall\; a\in A_1,\; n\in N$.
\ede
We can now say that a {\em strong connection on $H$\gal\ $B\inc P$
is a left $(B, (H^*)^{op})$-bimodule connection on $P$ vanishing on $1$}.
In an analogous way, we can define a right bimodule connection $\nabla_R$.
Then we can put them together and, in the spirit of \cite[p.284]{cq95},
define a bimodule connection as:
\bde\label{bicon}
Let $N$ be an $(A_1, A_2)$-bimodule and $\nabla_L$ and $\nabla_R$
a left and right bimodule connection, respectively.
We call a pair $(\nabla_L , \nabla_R)$ a {\em bimodule connection} on $N$.
\ede
Reasoning precisely as in \cite{cq95}, one can show that
an   $(A_1, A_2)$-bimodule $N$ admits a bimodule connection
if and only if it is projective as a bimodule
(i.e., as a module over $A_1\ot A_2^{op}$).
In a similar fashion, one can see that strong connections
correspond to equivariant connections in the algebraic-geometry setting
\cite[(20)]{r-d98}.

\bre\em\label{cgal}
Within the framework of the Hopf-Galois theory the right coaction
$\id\ot\dr : B\ot P\ra B\ot P\ot H$ and the restriction $\hD_{B\ot P}$
of the diagonal coaction $\hD_{P\ot P}$ (\ref{dsr}) coincide.
Therefore one can use either
of them to define the colinearity of a splitting $s$ of the multiplication map
$B\ot P\ra P$. In the general setting of $C$\gal s \cite[Definition 2.3]{bh99},
the diagonal
coaction $\hD_{P\ot P}:P\ot P\ra P\ot P\ot  C$ (coinciding with (\ref{dsr})
for \hge s) can be defined by the formula
$
\hD_{P\ot P}=(\id\ot\psi)\ci(\dr\ot\id)
$
\cite[Proposition 2.2]{bm98a},
where $\psi:C\ot P\ra P\ot C$ is an entwining structure and $\dr:P\ra P\ot C$,
$\dr(p)=p\0\ot p\1$, a coaction (see \cite[Section 3]{bm} for details).
If $B$ is the subalgebra of $P$ of $C$-coinvariants, i.e.,
$B=\{ b\in P\;|\;\dr(bp)=b\dr(p),\;\forall\; p\in P\}$, 
then \beq \hD_{P\ot P}(b\ot p) =(\id\ot\psi)(\dr(b1)\ot p)
=(\id\ot\psi)(b1\0\ot 1\1\ot p) =b1\0\ot\psi(1\1\ot p). \eeq On
the other hand, if $B\inc P$ is $C$-Galois and $\psi$ is its
canonical entwining structure \cite[(2.5)]{bh99}, then, by
\cite[Theorem 2.7]{bh99}, $P$ is a $(P,C,\psi)$-module
\cite{b-t99}, so that we have $\dr(p'p)=p'\0\psi(p'\1\ot p)$. In
particular, $\dr(p)=1\0\psi(1\1\ot p)$. Hence \beq
(\id\ot\dr)(b\ot p)=b\ot 1\0\psi(1\1\ot p). \eeq
 Therefore we need to distinguish between
$\hD_{B\ot P}$ and $\id\ot\dr$ in the $C$-Galois case.
\footnote{
We are grateful to T.\ Brzezi\'nski for suggesting to us this way of arguing.
}

If we define a strong connection on a $C$\gal\ $B\inc P$ as a unital
left $B$-linear right $C$-colinear (with respect to $\hD_{B\ot P}$)
splitting of the multiplication map $B\ot P\ra P$, then such a strong connection
yields a connection in the sense of \cite[Definition 3.5]{bm}. Indeed, let
$s$ be such a splitting, and $\Pi^s(r\d p):=r(s(p)-p\ot 1)$. One can verify
that this formula gives a well-defined left $P$-linear endomorphism of $\hO^1P$.
Furthermore, by the left $B$-linearity of $s$, for any
$\sum_i\d b_i.p_i\in(\hO^1B)P$, we have:
\beq
\Pi^s(\sum_i\d b_i.p_i)
= \sum_i\Pi^s\llp\d (b_i.p_i)-b_i\d p_i\lrp
= \sum_is(b_ip_i)- b_ip_i\ot 1-b_i(s(p_i)- p_i\ot 1)=0.
\eeq
Hence $P(\hO^1B)P\inc \ker\Pi^s$ by the left $P$-linearity of $\Pi^s$.
On the other hand, since $m\ci s=\id$ and $s(P)\inc B\ot P$, we have
$\pi_B(s(p))=1\ot_Bp$, where $\pi_B:P\ot P\ra P\ot_BP$ is the canonical
surjection. Consequently,
\beq
\pi_B(\Pi^s(p'\d p))
=\pi_B(p'(s(p)-p\ot 1))
=p'\pi_B(s(p))-rp\ot_B 1
=p'\ot_Bp-rp\ot_B 1
=\pi_B(p'\d p).
\eeq
Therefore, since $P(\hO^1B)P= \ker\pi_B$ (see below (\ref{dia})),
we obtain $\ker\Pi^s\inc P(\hO^1B)P$.
Thus $\ker\Pi^s= P(\hO^1B)P$. Next, take any $p\in P$. It follows from
$s(P)\inc B\ot P$ that
\beq
\d p-\Pi^s(\d p)=1\ot p-p\ot 1-s(p)+p\ot 1=1\ot p-s(p)\in B\ot P.
\eeq
Since also $m(1\ot p-s(p))=0$, we have
$\d p-\Pi^s(\d p)\in (\hO^1B)P\inc \ker\Pi^s$.
By the left $P$-linearity of $\Pi^s$ we can conclude now that
$\Pi^s\ci(\id-\Pi^s)=0$, i.e., $(\Pi^s)^2=\Pi^s$.
It remains to show that
$\hD_{P\ot P}\ci\Pi^s\ci\d=\llp(\Pi^s\ci\d)\ot\id\lrp\ci\dr$.
The property $\psi(c\ot 1)=1\ot c$ entails
\beq
\hD_{P\ot P}(p\ot 1)
= p\0\ot\psi(p\1\ot 1)
= p\0\ot 1\ot p\1.
\eeq
Therefore
\bea
\hD_{P\ot P}(\Pi^s(\d p))
\!\!\!\!\!\!&&
=\hD_{P\ot P}(s(p))-\hD_{P\ot P}(p\ot 1)
\nonumber\\ &&
=s(p\0)\ot p\1-p\0\ot 1\ot p\1
\nonumber\\ &&
=\llp(\Pi^s\ci\d)\ot\id\lrp(\dr(p))
\eea
by the colinearity of $s$. Consequently $\Pi^s$ is a connection, as claimed.
\ere

To exemplify Proposition~\ref{claim0} and Theorem~\ref{claim1},
let us translate the strong and non-strong connection forms
on quantum projective space $\IR P_q^2$ \cite[Example 2.8]{h-pm96}
to the language of splittings.
\bex\label{rp2}{\it (Quantum projective space $\IR P^2_q$.)~}\em
First let us recall how to define the coordinate ring $A (S_{q,\infty}^2)$
of the equator Podle\'s quantum sphere \cite{p-p87}.
To this end, we modify the convention in \cite{h-pm96}
by replacing $q$ by $q^{-1}$ and rewriting the generators as follows:
\beq
x=x_{11}, ~~~ y=x_{12}, ~~~ z= \frac{\sqrt{2(1+q^4)}}{1+q^2} x_{13}\; .
\eeq
Now we can define
$A(S_{q,\infty}^2)$ as $\IC \< x,y,z\> / I_{q,\infty}$,
where $\IC \< x,y,z\>$ is the (unital) free algebra generated by $x,y,z$ and
$I_{q,\infty}$ is the two-sided ideal generated by
\beq
x^2+y^2+z^2-1, \; xy-yx -i \mbox{$\frac{q^4 -1}{q^4 +1}$}z^2,\;
xz-\mbox{$\frac{q^2 +q^{-2}}{2}$}zx -i \mbox{$\frac{q^{-2} -q^{2}}{2}$}zy,\;
yz-\mbox{$\frac{q^2 +q^{-2}}{2}$}zy -i \mbox{$\frac{q^{2} -q^{-2}}{2}$}zx\; .
\eeq
To make $A(S_{q,\infty}^2)$ into a $\map (\IZ_2 , \IC )$-comodule algebra
we use the formulas (see above Section~6 in~\cite{p-p87} for the related
quantum-sphere automorphisms)
\beq
\dr (x) = x\ot \hg,\; \dr (y) = y\ot \hg,\; \dr (z) = z\ot \hg,
\eeq
where $\hg (\pm 1) = \pm 1$.
The coordinate ring of quantum projective space $\IR P_q^2$ is then defined
as the $\map (\IZ_2 , \IC )$-coinvariant subalgebra of  $A(S_{q,\infty}^2)$.
(The algebra $A(\IR P_q^2)$ is the subalgebra of $A(S_{q,\infty}^2)$
generated by the monomials of even degree.)
The extension $A(\IR P_q^2)\inc A(S_{q,\infty}^2)$ is a
$\map (\IZ_2 , \IC )$\gal\ which is  not cleft.
The non-cleftness can be proved by reasoning exactly as in \cite[Appendix]{hm99}.
Indeed, since $1$ and $\hg$ are linearly independent group-likes
($\hD c = c\ot c$), a cleaving map
$\Phi : \map (\IZ_2 , \IC )\ra  A(S_{q,\infty}^2)$
would have to map them to linearly independent (injectivity of $\Phi$)
invertible (convolution invertibility of $\Phi$) elements in
$A(S_{q,\infty}^2)$.
But this is impossible as $A(S_{q,\infty}^2)\inc \slq$,
and the only invertible elements in $\slq$ are non-zero numbers
\cite[Appendix]{hm99}.
Translating the formula \cite[Proposition 2.14]{h-pm96} for a strong connection
to our setting, we have
\beq
\ho (\hg ) = x\d x+y\d y +z\d z = x\ot x + y\ot y + z\ot z -1\ot 1 .
\eeq
The splitting corresponding to $\ho$ is then, due to its unitality
and left $A(\IR P_q^2)$-linearity, determined by
\beq
s(x)= x t~, ~~~s(y)= y t~, ~~~s(z)= z t~, ~~~\mbox{where }~~~
t = x\ot x + y\ot y + z\ot z \ .
\eeq
Thus one can directly see that the image of $s$ is in
$A(\IR P_q^2)\ot A(S_{q,\infty}^2)$.

Next,  consider a non-strong connection
$
\widetilde{\ho}(\hg ) =
\ho (\hg ) - 2\d x^2
$
\cite[Proposition 2.15]{h-pm96}.
Again, we compute the corresponding splitting:
\beq
\widetilde{s}(x)= s(x)-2x\d x^2~, ~~
\widetilde{s}(y)= s(y)-2y\d x^2~, ~~
\widetilde{s}(z)= s(z)-2z\d x^2\ .
\eeq
As in the proof of
\cite[Proposition 2.15]{h-pm96},
we can invoke the representation theory contained in 
\cite{p-p87} to conclude that $x\d x^2\neq 0$.
Consequently,
\beq
\llp(\id\ot \mbox{flip})\ci(\dr\ot \id -\id\ot\id\ot\id)\lrp (x\d x^2)=
x\d x^2\ot (\hg -1) \neq 0 .
\eeq
Hence the image of $\widetilde{s}$ is not in
$A(\IR P_q^2)\ot A(S_{q,\infty}^2)$.
\eex
\bre\em
Let $x,y,z$ be as above.
Since $x^2 + y^2 + z^2 =1$ (which was the reason for rescaling
the generators) and, \wrt\ the star structure
inherited from $SU_q(2)$, we have $ x^* = x, ~~y^*=y, ~~z^*=z $,
we can treat the generators $x,y,z$
as the Cartesian coordinates of $S_{q,\infty}^2$.
Having this in mind, we take the idempotent
$F = (x, y, z)^T (x, y, z) \in M_3(A(S^2_{q,\infty}))$
(here $^T$ stands for the matrix transpose) and define the projective module
of the normal bundle of $S^2_{q,\infty}$ as $A(S^2_{q,\infty})^3F$.
Therefore one can define the projective module of the
{\em tangent bundle of the  equator Podle\'s quantum sphere}
as $ A(S_{q,\infty}^2)^3  (\mbox{I}_3 - F)$, where $I_3$ is the identity matrix
in $M_3(A(S^2_{q,\infty}))$.
\ere

Let us now consider strong connections on {\em principal homogeneous
\hge s}, i.e., $P/I$\gal s given by a Hopf ideal $I$ in a Hopf algebra $P$.
Here the coaction is given by the formula
$\dr = (\id\ot\pi_I ) \ci \hD $, where $\pi_I$ is the canonical surjection
$P\ra P/I $.
For such extensions, it is known (e.g., see \cite[Theorem~2.1]{dhs99})
that if $B = P ^{co~P/I}$ then $I= B^+ P$, where $B^+ =\ker \he \cap B$.
If $s$ is the splitting associated to a strong connection, then,
due to the left $B$-linearity of $s$,
\beq\label{desc}
s( B^+ P ) = B^+ s(P) \inc B^+ B \ot P = B^+\ot  P\ .
\eeq
Hence $s$ descends to a splitting $i$ of the canonical surjection
$P \ra P/(B^+ P)$:
\vspace*{3mm}\beq\label{}
\def\normalbaselines{\baselineskip30pt\lineskip3pt\lineskiplimit3pt}
\def\mapright#1{\smash{\mathop{-\!\!\!\lra}\limits^{#1}}}
\def\Mapright#1{\smash{\mathop{-\!\!\!-\!\!\!-\!\!\!-\!\!\!\lra}\limits^{#1}}}
\def\mapleft#1{\smash{\mathop{\longleftarrow\!\!\!-}\limits_{#1}}}
\def\mapdown#1{\Big\downarrow\rlap{$\vcenter{\hbox{$\scriptstyle#1$}}$}}
\def\mapup#1{\Big\uparrow\llap{$\vcenter{\hbox{$\scriptstyle#1$}}$}}
\matrix{
&&P & \stackrel{\large\Mapright{s}}{\mapleft{m}} & B\ot P &&\cr
&& \mapdown{} && \mapdown{} && \cr
H&=&P/(B^+ P) & \stackrel{\large\Mapright{i}}{\mapleft{}}
& (B\ot P)/(B^+ \ot P) &=&P\ . \cr
}
\eeq
Explicitly, we have $i(\overline{p}) = \llp (\he\ot\id )\ci s \lrp (p)$.
(The map is well-defined because of (\ref{desc}).)
Put $s(p) = s(p)\z \ot s(p)\o$ (summation understood).
Then, as $m\ci s = \id$ and, for $b\in B$, $p\in P$,
$\he (b) p =bp ~\mbox{mod}~B^+ P$, we have
\beq
(\pi_I \ci i) (\overline{p}) = \pi_I \llp \he (s(p)\z ) s(p)\o \lrp =
 \pi_I \llp s(p)\z  s(p)\o \lrp = ( \pi_I \ci m \ci s ) (p) =\overline{p}\ .
\eeq
Furthermore, since $s$ is unital, so is $i$.
The right colinearity of $i$ follows from the strongness
($s(P) \inc B\ot P$) and the right colinearity of $s$:
\bea
(\dr \ci i) (\overline{p})
\!\!\!\!\!\!&&
= \he \llp s(p)\z \lrp {s(p)\o}\0 \ot {s(p)\o}\1
\nonumber\\ &&
= \llp (\he\ot\id\ot\id)\ci \hD_{P\ot P} \ci s \lrp (p)
\nonumber\\ &&
= \llp (\he\ot\id)\ci s \lrp (p\1)\ot \overline{p\2 }
\nonumber\\ &&
= i ( \overline{p\1 })\ot \overline{p\2 }
= i ( \overline{p}\1 )\ot \overline{p}\2 \ .
\eea
Thus one can associate to any strong connection on a principal homogeneous
\hge\ a total
integral of Doi \cite{d-y85} (unital right colinear map $H\ra P$).
Recall that total integrals always exist on faithfully flat \hge s
(\cite[Theorem~1]{s-hj90a}, \cite[(1.6)]{d-y85}, \cite[Remark~3.3]{s-hj90a}).
This is in agreement with point 3 of Corollary~\ref{ff}, although we claim
there only the left faithful flatness, and faithfully flat \hge s $B\inc P$
are defined as \hge s such that $P$ is $B$-faithfully-flat on both sides.
Note also that we could equally well proceed as in
\cite[Proposition~3.6]{bm98b} and define $i$ via a connection form.
If $s$ is the splitting associated to a connection form $\ho$,
i.e., $s=J_4 (\ho )$ (see (\ref{sf})), then
\bea
i (\overline{p})
\!\!\!\!\!\!&&
= \llp (\he \ot\id )\ci J_4 (\ho )\lrp (p)
\nonumber\\ &&
= (\he \ot\id)\llp p\ot 1 + p\1 \ho (\overline{p\2} ) \lrp
\nonumber\\ &&
=\he (p) \ot 1 + \he (p\1 ) \he (\ho (\overline{p\2 })^{(1)})\ot
\ho (\overline{p\2 })^{(2)}
\nonumber\\ &&
= \he_H (\overline{p}) \ot 1 + \he (p\1 ) (\he\ot\id ) (\ho (\overline{p\2 }))
\nonumber\\ &&
= \he_H (\overline{p}) \ot 1 + \llp (\he\ot\id )\ci \ho \lrp (\overline{p})\ ,
\eea
where $\ho (h) = \ho (h)^{(1)}\ot \ho (h)^{(2)}$,
summation understood, and $\he_H$ denotes the counit on $H$.
(See \cite[Proposition~3.6]{bm98b} for this kind of splittings in the case
of non-universal calculus.)
If $i$ is also left colinear, then, by \cite[Proposition~2.4]{hm99},
the formula $\ho = (S * \d )\ci i $ associates to $i$ a strong connection.
(Such connections are called canonical strong connections.)
It turns out that applying the above described way of associating a total
integral to a strong connection in the canonical case
is simply solving the equation $\ho = (S * \d )\ci i $ for $i$.
Indeed, since $\ho (h) = S i(h)\1 \d i(h )\2 $, we have
\beq\label{j4}
J_4(\ho )(p)= p\ot 1 + p\1 \ho (\overline{p\2} ) =
p\1 S i (\overline{p\2} )\1 \ot i(\overline{p\2} )\2\, .
\eeq
Applying $\he \ot \id$ yields
\beq\label{sol}
i (\overline{p} ) = \llp (\he \ot \id )\ci J_4(\ho ) \lrp (p) \ ,
\eeq
as claimed.

\bex\label{q}{\it (Quantum and classical Hopf fibration.)~}\em
The above described formalism applies to the quantum Hopf fibration.
We refer to \cite{hm99} for the computation of
projector matrices of the quantum Hopf line bundles
from the Dirac q-monopole connection \cite{bm93},
and to \cite{h-pm} for the computation of the Chern-Connes pairing
of these matrices with the cyclic cocycle (trace) \cite[(4.4)]{mnw91}.
(This pairing yields numbers called ``Chern numbers"" or ``charges."
See Proposition~\ref{free} for the freeness of the direct sum of charge $-1$
and charge $1$ quantum Hopf line bundles,
cf.\ \cite[(4.2)]{ds94} for a local description of such bundles.)
Here we only remark that this quantum principal fibration
admits infinitely many canonical strong connections.
Indeed, for the injective antipode (which is the case here),
\cite[Corrolary 2.6]{hm99} classifies
the canonical strong connections by unital bicolinear splittings.
On the other hand, by \cite[p.363]{mmnnu91}, all unital bicolinear splittings
$i : \IC [z, z^{-1}] \ra \slq$ are of the form:
\bea
\label{split}
&&
i (z^n) = (1+ \hz p_n(\hz)) \ha^n ,
\nonumber \\ &&
i (z^{-n}) = (1+ \hz r_n(\hz)) \hd^n ,
\eea
where $p_n, r_n$ are arbitrary polynomials in $\hz : = -q^{-1} \hb \hg$.
(Here, $\ha, \hb, \hg, \hd$ are the generators of \slq\ as in \cite{hm99}.)
Since the equality $\ho = (S * \d ) \ci i$ can be solved for $i$
(see (\ref{sol})), different splittings yield different connections.
Hence there are infinitely many connections.

However, for $q=1$, after passing to the
de Rham forms, all the canonical strong connections
coincide with the classical Dirac monopole.
More precisely, let $\pi_{DR}$ be the canonical
projection from the universal onto the de Rham differential calculus and
$i_0$ be the splitting corresponding to the Dirac monopole
(i.e., given by (\ref{split}) with $p_n=0=r_n$ for all $n$).
Then
\beq
\pi_{DR}\ci (S * \d ) \ci i = \pi_{DR}\ci (S * \d ) \ci i_0 ~~\mbox{for all}~~ i ~.
\eeq
Indeed, we have
\beq
\llp (S*\d_{DR})\ci (i-i_0 )\lrp(z) =  (S*\d_{DR})(\hb\hg\, p_1\! (\hb\hg)\ha) \ .
\eeq
Furthermore, using the commutativity of functions with forms and functions,
and the Leibniz rule, we obtain
\bea\label{drr}
 (S*\d_{DR})(hh')
\!\!\!\!\!\!&&
=S(h\1)S(h'\1) \d_{DR}(h\2h'\2)
\nonumber \\ &&
=S(h'\1)h'\2 S(h\1) \d_{DR}h\2 + S(h\1)h\2 S(h'\1) \d_{DR}h'\2
\nonumber \\ &&
=\he(h') S(h\1) \d_{DR}h\2 + \he(h) S(h'\1) \d_{DR}h'\2 \ .
\eea
Substituting $h=\hb$ and $h' = \hg\, p_1\! (\hb\hg)\ha$, and noting that
$ \he(\hb)=0$ and $\he (\hb\hg\, p_1\! (\hb\hg)\ha)= 0$,
one can conclude that
$(S * \d_{DR} )(i(z)) = (S * \d_{DR} )(i_0 (z)) $.
On the other hand, for any connection form $\ho$ we have
\beq
(\pi_{DR}\ci\ho ) (uu') = (\pi_{DR}\ci\ho )(u)\he(u') +
\he(u) (\pi_{DR}\ci\ho )(u')\ .
\eeq
Therefore $(S * \d_{DR} )\ci i$ and $(S * \d_{DR} )\ci i_0$ coincide
on any power of $z$, whence are equal, as claimed.
\eex

\section{Chern-Connes pairing for the super Hopf fibration}
\setcounter{equation}{0}

The super Hopf fibration leading to the super sphere has an interesting history.
To the best of our knowledge, it was first introduced by Landi and Marmo
\cite{lm87}. They treated supersymmetric abelian gauge fields in general
and worked out details for the super group $UOSP(1,2)$.
Everything was formulated within the Grassmann envelope of the super algebra
$uosp(1,2)$. The super manifold theory has been used in the work
of Teofilatto \cite{t-p88}.
He defines and studies super Riemann surfaces.
As the simplest example, he treated the super sphere with $S^2$ as its body.
Ideas of noncommutative geometry were used in  \cite{gkp96,gkp97}
to introduce an ultraviolet regularization for quantum fields defined on $S^2$.
The fuzzy sphere \cite{m-j92} was introduced through suitable embeddings
of the algebra of $N\times N$ matrices. In \cite{gkp96},
similar embeddings of modules led to approximation of sections of line bundles
over $S^2$. Also in \cite{gkp96}, there is a study of fermions
and supersymmetric extensions of the fuzzy sphere.
An extensive treatment of the approximation of super-graded functions over
the super sphere, and sections of a bundle through sequences of graded modules,
as well as the treatment of the graded de Rham complex, is given in \cite{gr98}.
The description of the monopole on the super sphere that we provide
can be related to that given in \cite{bbl90}.
A detailed study of the super monopole using the super-geometry approach
can be found in \cite{l-gb}.

Our approach here to the super Hopf fibration is purely algebraic.
First, we show that the super Hopf fibration can be considered as
an $H$-Galois extension $\bs\inc\as$, where $H = \IC [z, z^{-1}]$
is the Hopf algebra generated by invertible group-like element $z$.
The polynomial algebras \as\ and \bs\
are taken as nilpotent extensions (by two Grassmann variables $\hl_{\pm}$)
of the (complex) coordinate rings
of 3-dimensional sphere  $S^3$ and 2-dimensional sphere  $S^2$, respectively
(see \cite{gkp96}).
This is summed up in the following commutative diagram with exact columns
(but not rows):

\beq\label{diag}
\def\normalbaselines{\baselineskip30pt
\lineskip3pt \lineskiplimit3pt }
\def\mapright#1{\smash{
\mathop{\!\!\!-\!\!\!\longrightarrow\!\!\!}
\limits^{#1}}}
\def\mapdown#1{\Big\downarrow
\rlap{$\vcenter{\hbox{$\scriptstyle#1$}}$}}
\matrix{
\bs\cap\< \hl_{\pm}\> & \mapright{} & \< \hl_{\pm}\> & & \cr
\mapdown{} && \mapdown{} &&  \cr
\bs & \mapright{} & \as & \mapright{} & H \cr
\mapdown{\wp } && \mapdown{} && \bigg\| \cr
A(S^2) & \mapright{} & A(S^3) & \mapright{} & H \cr
}
\eeq

\noindent
Thus, in a sense, the  super Hopf fibration can be viewed as a Grassmann
covering of the classical (complex) Hopf fibration.

\bde\label{shdef}
Let $R = \IC [a, b, c, d]$ be the polynomial ring in four variables.
Put $D = ad - bc$.
Let $I$ be the two-sided ideal in the (unital) free algebra
$R\langle \hl_+ , \hl_- \rangle $ generated by
\beq
\hl_{+}^2 ,~~ \hl_{-}^2 ,~~ \hl_+ \hl_- + \hl_- \hl_+ ,
~~\hl_+ \hl_- + D - 1 .
\eeq
We call the quotient algebra $A(S^3_s) := R\langle \hl_+ , \hl_- \rangle / I$
the coordinate ring of 3-dimensional super sphere  $S^3_s$.
\ede
It can be easily verified that the (matrix) formula
\beq
\dr \pmatrix{a & b \cr c & d \cr \hl_+ & \hl_-\cr}
=
\pmatrix{a\ot1 & b\ot1 \cr c\ot1 & d\ot1 \cr
\hl_+\ot1 & \hl_-\ot1\cr}
\pmatrix{1\ot z & 0 \cr 0 & 1\ot z^{-1} \cr}
\eeq
defines a coaction
$\dr : \as \ra \as\ot H$
making \as\ a right $H$-comodule algebra.
\ble\label{B}
Let $\bs := \{a\in \as ~|~ \dr (a) = a \ot 1\}$
be the algebra of $H$-coinvariants.
Then
\bs\ is the subalgebra of \as\ generated by
\beq\label{gen}
1,~~ ab,~~ bc,~~cd, ~~\hl_{+}b,~~\hl_{+}d, ~~\hl_{-}a,~~\hl_{-}c,
~~ \hl_+ \hl_- .
\eeq
\ele
\bpf
Evidently, the algebra generated by (\ref{gen}) is contained in \bs.
For the opposite inclusion, note first that every non-zero
element $a$ of \as\ can be written as a linear combination of
non-zero monomials $m_{k,\ell}$ such that $\dr (m_{k,\ell}) = m_{k,\ell} \ot z^k$.
Since the powers of $z$ form a basis of $H$,
if $a\in \bs$, then
$a$ must be a linear combination of
non-zero monomials $m_{0,\ell}$.
On the other hand, any $m_{0,\ell}$ is a word composed
of the same number of letters coming from the alphabet
$\{a, c, \hl_+\}$ and the alphabet $\{b, d, \hl_-\}$.
Furthermore, since all letters commute or anti-commute,
we can always pair the letters coming from different alphabets.
Hence $m_{0,\ell}$ can be expressed in terms of (\ref{gen}),
as needed.
\epf
\vspace*{-1cm}
\bpr\label{hgs}
The extension of algebras
$\bs \inc \as $ is $H$-Galois.
\epr
\bpf
Define the map $\tau : H \ra \as \ot_{A(S^2_s)} \as$
by the formulas ($n\in \IN$):
\bea\label{tau}
&&
\tau (z^n) = (1+n\hl_+ \hl_- ) \sum_{k=0}^{n}{\scriptsize\pmatrix{n\cr k\cr}}
d^{n-k} (-b)^k  \ot_{A(S^2_s)} a^{n-k} c^k , \nonumber \\
&&
\tau (z^{-n}) = (1+n\hl_+ \hl_- ) \sum_{k=0}^{n}{\scriptsize\pmatrix{n\cr k\cr}}
a^{n-k} (-c)^k  \ot_{A(S^2_s)} d^{n-k} b^k.
\eea
We are going to prove that
$\tilde\chi := (m\ot \id) \ci (\id \ot \tau )$
is the inverse of the canonical map $\chi$.
(This means that $\tau $ is the translation map.)
Since $\chi $ and $\tilde\chi$ are both left
\as-linear maps by construction,
it suffices to check $\chi \ci \tilde\chi= \id$
and  $\tilde\chi \ci \chi= \id$ on elements
of the form $1\ot h$ and $1\ot_{A(S^2_s)} p$, respectively.
To verify the first identity, we recall that any $h \in H$
is a linear combination of $z^{\pm n}, \; n\in \IN$,
and compute:
\bea
(\chi \ci \tilde\chi)(1\ot z^n)
\!\!\!\!\!\!&&
=(1 + n\hl_+ \hl_- ) \mbox{$\sum_{k=0}^n{\scriptsize\pmatrix{n\cr k\cr}}
d^{n-k} (-b)^k \chi(1 \ot_{A(S^2_s)} a^{n-k} c^k )$}\nonumber \\
&&
=(1 + n\hl_+ \hl_- )\llp\mbox{$\sum_{k=0}^{n}{\scriptsize\pmatrix{n\cr k\cr}}
d^{n-k} (-b)^k  a^{n-k} c^k$}\lrp\ot  z^n  \nonumber \\
&&
=(1 + n\hl_+ \hl_- )(ad-bc)^n \ot z^n  \nonumber \\
&&
=(1 + n\hl_+ \hl_- )(1 - \hl_+ \hl_-)^n \ot z^n  \nonumber \\
&&
=(1 + n\hl_+ \hl_- )(1 - n\hl_+ \hl_-)\ot z^n  \nonumber \\
&&
=1\ot z^n.
\eea
In the fourth equality we used the determinant relation
$\hl_+ \hl_- + ad -bc = 1$. Similarly, we obtain:
\beq
(\chi \ci \tilde\chi) (1\ot z^{-n})
 =(1 + n\hl_+ \hl_- ) \mbox{$\sum_{k=0}^n{\scriptsize\pmatrix{n\cr k\cr}}
a^{n-k} (-c)^k \chi(1 \ot_{A(S^2_s)} d^{n-k} b^k )$}
=1\ot z^{-n}.
\eeq
Thus $\chi \ci \tilde\chi = \id$.
For the other identity, we note that it is sufficient to check it on the monomials
$m_{\pm n}$ (as in the proof of Lemma~\ref{B} but with the second index
suppressed). Since
$\dr (m_{\pm n}) = m_{\pm n} \ot z^{\pm n}$, we have
$m_{n} d^{n-k} b^k \in \bs$ and $m_{-n} a^{n-k} c^k \in \bs$.
Hence, using the centrality of $\hl_+ \hl_- \in\bs $, we can compute:
\bea
(\tilde\chi \ci \chi)(1\ot_{A(S^2_s)} m_n)
\!\!\!\!\!\!&&
= m_n \tilde\chi (1 \ot z^n )
\nonumber \\ &&
=m_n (1 + n\hl_+ \hl_- ) \mbox{$\sum_{k=0}^n{\scriptsize\pmatrix{n\cr k\cr}}
d^{n-k} (-b)^k \ot_{A(S^2_s)} a^{n-k} c^k $}\nonumber \\
&&
=(1 + n\hl_+ \hl_- )\mbox{$\sum_{k=0}^{n}{\scriptsize\pmatrix{n\cr k\cr}}
m_n d^{n-k} (-b)^k \ot_{A(S^2_s)} a^{n-k} c^k$}    \nonumber \\
&&
= 1 \ot_{A(S^2_s)}(1 + n\hl_+ \hl_- )
\mbox{$\sum_{k=0}^{n}{\scriptsize\pmatrix{n\cr k\cr}}
m_n d^{n-k} (-b)^k a^{n-k} c^k$}    \nonumber \\
&&
= 1 \ot_{A(S^2_s)} m_n (1 + n\hl_+ \hl_- )(ad-bc)^n \nonumber \\
&&
=1\ot_{A(S^2_s)} m_n.
\eea
Here the last step is as in the previous calculation.
Similarly, we get:
\bea
(\tilde\chi \ci \chi)(1\ot_{A(S^2_s)} m_{-n})
\!\!\!\!\!\! && 
=(1 + n\hl_+ \hl_- )\mbox{$\sum_{k=0}^{n}{\scriptsize\pmatrix{n\cr k\cr}}
m_{-n} a^{n-k} (-c)^k \ot_{A(S^2_s)} d^{n-k} b^k$}
\nonumber\\ &&
= 1\ot_{A(S^2_s)} m_{-n}.
\eea
Therefore
$\tilde\chi$ is the inverse of $\chi$, and the extension is $H$-Galois.
\epf
We use the idea of colinear lifting (\ref{dia3}) to construct a connection form.
We consider this connection as the (universal-calculus) super Dirac monopole.
Since it is strong, we can conclude that
the extension $\bs \inc \as $ enjoys all properties itemized in
Corollary~\ref{ff}.
\bpr\label{omprop}
Let  $\ho : H \ra \hO^1 \as$ be the linear map defined by ($n\in \IN$)
\bea
&&\label{omform-}
\ho (z^n)
= (1+n\hl_+ \hl_- ) \sum_{k=0}^{n}{\scriptsize\pmatrix{n\cr k\cr}}
d^{n-k} (-b)^k  \d (a^{n-k} c^k) , \\
&&\label{omform+}
\ho (z^{-n}) = (1+n\hl_+ \hl_- ) \sum_{k=0}^{n}{\scriptsize\pmatrix{n\cr k\cr}}
a^{n-k} (-c)^k \d ( d^{n-k} b^k) \ .
\eea
Then \ho\ is a strong connection form.
\epr
\bpf
Note first that $\ho (1) = 0$
and $\dsr \ho (z^{\pm n}) =\ho (z^{\pm n})\ot 1$.
Furthermore,
\beq
\llp (m \ot \id )\ci (\id \ot \dr )\ci\ho \lrp (z^{\pm n})
=\chi \llp \tau (z^{\pm n}) - 1\ot_{A(S^2_s)} 1\lrp
=1\ot \llp z^{\pm n} - \he (z^{\pm n})\lrp .
\eeq
This proves that \ho\ is a connection form.
It remains to check the strongness condition.
By linearity, it suffices to do it on monomials $m_{\pm n}$
(see the proof of Proposition~\ref{hgs}).
Putting
\bea
&&
p_k^+=(1+n\hl_+ \hl_- ){\scriptsize\pmatrix{n\cr k\cr}} d^{n-k} (-b)^k,
~~~ q_k^+= a^{n-k} c^k,
\nonumber \\ &&
p_k^-=(1+n\hl_+ \hl_- ){\scriptsize\pmatrix{n\cr k\cr}} a^{n-k} (-c)^k,
~~~ q_k^-= d^{n-k} b^k ,
\eea
 and using the Leibniz rule, we obtain:
\bea
&&
\phantom{=}\d m_{\pm n} - m_{\pm n}\ho (z^{\pm n})
\nonumber \\ &&
= \d m_{\pm n} - \mbox{$\sum_{k=0}^{n} m_{\pm n}$} p_k^\pm \d q_k^\pm
\nonumber \\ &&
= \d m_{\pm n} - \d \llp m_{\pm n} \mbox{$\sum_{k=0}^{n} $} p_k^\pm  q_k^\pm \lrp
+\mbox{$\sum_{k=0}^{n}$} \d ( m_{\pm n}  p_k^\pm ).  q_k^\pm
\nonumber \\ &&
= \mbox{$\sum_{k=0}^{n}$} \d ( m_{\pm n}  p_k^\pm ).  q_k^\pm
 ~\in~ (\hO^1  \bs )\as.
\eea
Consequently, for any $a\in \as$, we have
$( \id - \Pi^{\omega} )(\d a)\in (\hO^1  \bs )\as $,
i.e., \ho\ is strong.
\epf
\vspace*{-1cm}

Our next step is to consider super Hopf line bundles.
Precisely as in the case of quantum Hopf line bundles
\cite[Definition~3.1]{hm99},
since we are dealing with one-dimensional corepresentations of $H$
($\rho_\mu (1) = 1\ot z^{-\mu}$), we can identify the colinear maps
$\xi : \IC\ra\as $ with their values at $1$
($\eta (\xi):= \xi (1)$), and define them as the following
bimodules over \bs:
\beq\label{amu}
\as_\mu := \{a\in \as ~|~ \dr a=a\ot z^{-\mu}\},~\mu \in\IZ .
\eeq
Reasoning in a similar manner as in the proof of Lemma~\ref{B},
one can see that ($n>0$)
\bea
&&\label{as-}
\as_{-n}
= \mbox{$\sum_{k=0}^{n}$}\bs a^{n-k} c^k +
 \mbox{$\sum_{k=0}^{n-1}$}\bs a^{n-1-k} c^k \hl_+ ,\\
&&\label{as+}
\as_{n}
= \mbox{$\sum_{k=0}^{n}$}\bs d^{n-k} b^k +
 \mbox{$\sum_{k=0}^{n-1}$}\bs d^{n-1-k} b^k \hl_- .
\eea
Note that, since the powers of $z$ form a basis of $H$,
we have the direct sum decomposition
$\as = \bigoplus_{\mu\in\IZ } \as_\mu $ as \bs-bimodules.
Observe also that the bimodules $\as_\mu$ provide examples of bimodules $P_\rho$
defined in Proposition~\ref{pero} (cf.\ \cite[Appendix B]{d-m97a}).
Our goal is to compute projector matrices of these modules
and their pairing with the appropriate cyclic cocycle on~\bs.
The strategy for computing the projector matrices is to use the
splitting associated to super Dirac monopole (Proposition~\ref{omprop})
and Lemma~\ref{genl}.
To apply the aforementioned lemma, first we need to show that
the monomials occurring in  formula (\ref{as-}) are linearly independent,
and that the same holds for the monomials in (\ref{as+}).
\ble\label{li}
\bea
&&\label{li-}
\sum_{k=0}^{n}\ha_k a^{n-k} c^k +
\sum_{\ell =0}^{n-1}\hb_\ell a^{n-1-\ell} c^\ell \hl_+ = 0 ~~~
\Rightarrow~~~ \ha_k =0=\hb_\ell ,
~~~\fa k, \ell \ ; \\
&&\label{li+}
\sum_{k=0}^{n}\ha_k d^{n-k} b^k +
\sum_{\ell =0}^{n-1}\hb_\ell d^{n-1-\ell} b^\ell \hl_- = 0 ~~~
\Rightarrow~~~ \ha_k =0=\hb_\ell ,
~~~\fa k, \ell \ .
\eea
\ele
\bpf
Let $\tilde R = \IC [\ta, \tb, \tc, \td ]/\<\ta \td - \tb \tc -1 \>$
denote the coordinate ring of $SL(2,\IC )$
and $\IC [\hl ]/\<\hl^2 \>$ be the algebra of dual numbers.
We have the following homomorphism of algebras:
\bea
&
\pi : \as \lra \tilde R\ot  \IC [\hl ]/\<\hl^2 \> ,
&\nonumber\\
&\pi (a) = \ta \ot 1,~\pi (b) = \tb \ot 1,~ \pi (c) = \tc \ot 1,~\pi (d) = \td \ot 1,
~\pi (\hl_\pm ) = 1  \ot \hl .&
\eea
Applying $\pi$ to the first equality in (\ref{li-}) yields
\beq
\sum_{k=0}^{n}\ha_k \ta^{n-k}\tc^k\ot 1 +
\sum_{\ell =0}^{n-1}\hb_\ell \ta^{n-1-\ell}\tc^\ell\ot \hl
= 0\ .
\eeq
Since the monomials $\ta^{n-k}\tc^k$ are part of the PBW basis of $\tilde R$,
they are linearly independent. Hence
$\ha_k =0=\hb_\ell ,\fa k, \ell $, by the linear independence of $1$ and \hl.
The second implication can be proved in the same way.
\epf
Note now that the above described identification $\eta$
allows one to identify $s_{\rho_\mu} $ of (\ref{sro}) with the restriction of $s$
to $\as_\mu$ (see Proposition~\ref{pero}):
\beq
\as_\mu \ni \widetilde{\xi} \to \llp (\id\ot\eta )\ci s_{\rho_\mu} \ci \eta^{-1}
\lrp (\widetilde{\xi})\in \bs \ot \as_\mu \ ,
\eeq
\beq
\llp (\id \ot\eta )\ci s_{\rho_\mu} \ci \eta^{-1}\lrp (\widetilde{\xi})
=
\overline{\ell }\llp s_{\rho_\mu} ( \eta^{-1} (\widetilde{\xi}))\lrp (1)
= \llp s\ci \eta^{-1}(\widetilde{\xi})\lrp (1) = s (\widetilde{\xi}) \ .
\eeq
On the other hand, remembering the formula for the universal differential
and using again the fact that $(ad-bc)^n = 1-n\hl_+ \hl_-$,
we can write (\ref{omform-}) in the following form:
\beq
\ho (z^n) = (1+n\hl_+ \hl_- )\sum_{\ell=0}^{n} {\scriptsize\pmatrix{n\cr \ell\cr}}
d^{n-\ell} (-b)^\ell  \ot a^{n-\ell} c^\ell - 1\ot 1 .
\eeq
Substituting this to (\ref{sf}), we obtain
\bea\label{s-}
s(a^{n-k} c^k)
\!\!\!\!\!\!&&
= a^{n-k} c^k\ot 1 + a^{n-k} c^k \ho (z^n)
\nonumber \\ &&
= \mbox{$\sum_{\ell=0}^{n}$} a^{n-k} c^k
(1+n\hl_+ \hl_- ) {\scriptsize\pmatrix{n\cr \ell\cr}}
d^{n-\ell} (-b)^\ell  \ot a^{n-\ell} c^\ell ,
\nonumber \\ &&
~\nonumber \\
s(a^{n-1-k} c^k \hl_+)
\!\!\!\!\!\!&&
= a^{n-1-k} c^k \hl_+ \ot 1 + a^{n-1-k} c^k \hl_+ \ho (z^n)
\nonumber \\ &&
= \mbox{$\sum_{\ell=0}^{n}$} a^{n-1-k} c^k \hl_+
{\scriptsize\pmatrix{n\cr \ell\cr}}
d^{n-\ell} (-b)^\ell  \ot a^{n-\ell} c^\ell .
\eea
Similarly, substituting
\beq
\ho (z^{-n}) = (1+n\hl_+ \hl_- )\sum_{\ell=0}^{n} {\scriptsize\pmatrix{n\cr \ell\cr}}
a^{n-\ell} (-c)^\ell  \ot d^{n-\ell} b^\ell - 1\ot 1
\eeq
to (\ref{sf}), we get
\bea\label{s+}
&&
s(d^{n-k} b^k)
= \mbox{$\sum_{\ell=0}^{n}$} d^{n-k} b^k
(1+n\hl_+ \hl_- ) {\scriptsize\pmatrix{n\cr \ell\cr}}
a^{n-\ell} (-c)^\ell  \ot d^{n-\ell} b^\ell ,
\nonumber \\ &&
s(d^{n-1-k} b^k \hl_-)
= \mbox{$\sum_{\ell=0}^{n}$} d^{n-1-k} b^k \hl_-
{\scriptsize\pmatrix{n\cr \ell\cr}}
a^{n-\ell} (-c)^\ell  \ot d^{n-\ell} b^\ell .
\eea
Hence, by Lemma~\ref{li} and Lemma~\ref{genl},
we can conclude that $\as_{-n} = \bs^{2n+1} E_{-n}$ as
left \bs-modules,
where $E_{-n} = P_{-n} Q_{-n}^T$
(symbol $^T$ stands for the matrix transpose) with
\bea
&&
P_{-n}^T := (1+n\hl_+ \hl_- ) \llp a^n,\cdots , a^{n-k} c^k ,\cdots , c^n,
 a^{n-1} \hl_+ ,\cdots , a^{n-1-k} c^k \hl_+ ,\cdots ,  c^{n-1} \hl_+ \lrp ,
\nonumber \\ &&
Q_{-n}^T := \llp d^n,\cdots , {\scriptsize\pmatrix{n\cr \ell\cr}} d^{n-\ell} (-b)^\ell ,
\cdots , (-b)^n, 0 ,\cdots , 0 \lrp .
\eea
In an analogous manner, we infer that $\as_{n} = \bs^{2n+1} E_{n}$
as left \bs-modules, where $E_{n} = P_{n} Q_{n}^T$ with
\bea
&&
P_{n}^T := (1+n\hl_+ \hl_- ) \llp d^n,\cdots , d^{n-k} b^k ,\cdots , b^n,
 d^{n-1} \hl_- ,\cdots , d^{n-1-k} b^k \hl_- ,\cdots ,  b^{n-1} \hl_- \lrp ,
\nonumber \\ &&
Q_{n}^T
:= \llp a^n,\cdots , {\scriptsize\pmatrix{n\cr \ell\cr}} a^{n-\ell} (-c)^\ell ,
\cdots , (-c^n), 0 ,\cdots , 0 \lrp .
\eea
To show the non-freeness of the above projective modules,
we determine the Chern-Connes pairing between their classes
in $K_0(\bs )$ and the cyclic cocycle on \bs\ obtained by
the pull-back $\wp ^*$ (see (\ref{diag})) of the cyclic 2-cocycle $c_2$
on $A(S^2)$ given by the integration on $S^2$.
We have:
\beq\label{CC}
\< \wp^* (c_2), [E_{\pm n}] \> = \< c_2 , [\wp_* E_{\pm n}] \> = \pm n .
\eeq
Here the last equality follows from the fact that matrix
$( \wp_* E_{\pm n})_{i,j} := \wp\llp ( E_{\pm n})_{i,j}\lrp $
is a projector matrix of the classical Hopf line bundle
with the Chern number  equal to $\pm n$.
Furthermore, since every free module can be represented in $K_0$ by the identity
matrix, the $K_0$-class of any free \bs-module always vanishes.
(The Chern class of a trivial bundle is zero.)
Thus the left modules $\as_\mu$, $\mu\neq 0$, are not (stably) free
and are pairwise non-isomorphic.
Now, reasoning as in \cite[Section~4]{hm99}, we obtain:
\bco\label{noncl}
The $H$\gal\ $\bs\inc\as$  (super Hopf fibration) is {\em not} cleft.
\eco
Let us remark that projectors $E_{\pm n}$ are not hermitian
\wrt\ the involution
\beq\label{*}
a^* = d, ~~b^*=-c, ~~c^*=-b, ~~d^*=a, ~~\hl_{\pm}^*=-\hl_{\mp} .
\eeq
Nevertheless, one can slightly modify $E_{\pm n}$ to find hermitian projectors
$F_{\pm n}=F_{\pm n}^{\dag}$ such that the modules
$\bs^{2n+1} E_{\pm n}$ and $\bs^{2n+1} F_{\pm n}$ are isomorphic.
They are given by the formulas $F_{\pm n} = U_{\pm n} U_{\pm n}^{\dag}$,
where ($n>0$)
\bea
U_{-n}^T := &&\!\!\!\!\!\!\!\!
(1+\frac{n-1}{2}\hl_+ \hl_- ) \times
\nonumber \\ && \!\!\!\!\!\!\!\!
\llp a^n,\cdots , {\scriptsize\pmatrix{n\cr k\cr}}^\frac{1}{2} a^{n-k} c^k ,
\cdots , c^n, a^{n-1} \hl_+ ,\cdots ,
{\scriptsize\pmatrix{n-1\cr k\cr}}^\frac{1}{2} a^{n-1-k} c^k \hl_+ ,\cdots ,
c^{n-1} \hl_+ \lrp,
\nonumber \\ &&
\nonumber \\
U_{n}^T := &&\!\!\!\!\!\!\!\!
(1+\frac{n+1}{2}\hl_+ \hl_- )  \times
\nonumber \\ && \!\!\!\!\!\!\!\!
\llp d^n,\cdots , {\scriptsize\pmatrix{n\cr k\cr}}^\frac{1}{2} d^{n-k} b^k ,
\cdots , b^n, d^{n-1} \hl_- ,\cdots ,
{\scriptsize\pmatrix{n-1\cr k\cr}}^\frac{1}{2} d^{n-1-k} b^k \hl_- ,\cdots ,
b^{n-1} \hl_- \lrp.
\eea
The matrices $ F_{\pm n}$ are hermitian by construction.
To check that they are idempotent, we compute:
\bea\label{UU}
U_{-n}^{\dag} U_{-n} \!\!\!\!\!\! && =
(1+\frac{n-1}{2}\hl_+ \hl_- )^2
\mbox{$\sum_{k=0}^{n}$} (ad)^{n-k} (-bc)^k
- \hl_- \hl_+ \mbox{$\sum_{k=0}^{n-1}$} (ad)^{n-1-k} (-bc)^k
\nonumber \\ &&
= (1+(n-1)\hl_+ \hl_- ) (ad-bc)^n + \hl_+ \hl_- (ad-bc)^{n-1}
\nonumber \\ &&
= (1+(n-1)\hl_+ \hl_- ) (1-n\hl_+ \hl_- ) +  \hl_+ \hl_-  (1-(n-1)\hl_+ \hl_- )
\nonumber \\ &&
=1 \ .
\eea
In the same manner, we check $U_{n}^{\dag} U_{n} =1 $.
It remains to verify that the projective modules
$\bs^{2n+1} E_{\pm n}$ and $\bs^{2n+1} F_{\pm n}$ are isomorphic.
For this purpose, we use (\ref{LL}) and take as $L, \tilde L$, the matrices
$L_{\pm n} := U_{\pm n}Q_{\pm n}^T,
~~\tilde L_{\pm n} := P_{\pm n}V_{\pm n}^T \in M_{2n+1}(\bs)$, respectively.
A calculation similar to (\ref{UU}) shows that $Q_{\pm n}^T P_{\pm n} = 1$.
This together with (\ref{UU}) and $U_{n}^{\dag} U_{n} =1 $
implies that $L_{\pm n}$ and $\tilde L_{\pm n} $ satisfies (\ref{ll}).
(Note that $L_{\pm n} \tilde L_{\pm n} = E_{\pm n}$ and
 $ \tilde L_{\pm n}L_{\pm n} = F_{\pm n}$.)
Thus the modules $\bs^{2n+1}E_{\pm n} $ and $\bs^{2n+1}F_{\pm n} $
are isomorphic, as claimed. This hermitian presentation of the
projective modules $\as_{\pm n}, \, n>0$, agrees with \cite[(3.25)]{l-ga}
for the projectors of the classical Hopf line bundles,
and resembles the appropriate formulas obtained in \cite[Section~4.2]{l-gb}.
(The case $n=0$ is trivial.)

Finally, we want to show that
\bpr\label{free}
$\as_{-1} \oplus \as_1 = \bs^2 $ as left $\bs$-modules.
\epr
\bpf
We can infer from the preceding considerations that the matrix
$\diag (F_{-1}, F_1 )$ is a projector matrix of $\as_{-1} \oplus \as_1$.
First, it turns out technically convenient to conjugate $F_1$ by
\beq
M := \pmatrix{0 & 1 & 0\cr 1 & 0 & 0\cr 0 & 0 & 1}\ .
\eeq
Then $\tilde F_1 : = M F_1 M$ is evidently equivalent (i.e., giving an isomorphic
projective module) to $F_1$ (just take $L = \tilde L = M$ in (\ref{ll})),
whence $\diag (F_{-1}, F_1 )$ is equivalent to $\diag (F_{-1}, \tilde F_1 )$.
(Note that this way we have $F_{-1} \tilde F_1 = 0 = \tilde F_1 F_{-1} $.)
To prove the proposition we employ (\ref{LL}-\ref{ll}),
and put $F = \diag (F_{-1}, \tilde F_1 )$ and $ E = \diag (1, 1)$.
The point is to find $L, \tilde L$ satisfying (\ref{ll}).
Since
\beq
F_{-1} = (a, c, \hl_+ )^T (d, -b, -\hl_- )
~~\mbox{and}~~
\tilde F_{1} = (1+ 2\hl_+ \hl_- )(b, d, \hl_- )^T (-c, a, -\hl_+ ) \ ,
\eeq
we look for  $ \tilde L$  of the form
\beq
 \tilde L = \pmatrix{f_- \cr f_+ }~, ~~
f_- = \pmatrix{a\cr c \cr \hl_+ }  \pmatrix{u_+ & v_+ }~, ~~
f_+ = \pmatrix{b\cr d \cr \hl_- }  \pmatrix{u_- & v_- }~,
\eeq
and for $L$  of the form
\beq
 L = \pmatrix{g_- & g_+ }~, ~~
g_- = \pmatrix{x_- \cr y_- } \pmatrix{d& -b& -\hl_- } ~, ~~
g_+ = \pmatrix{x_+ \cr y_+}  \pmatrix{-c & a& -\hl_+ }~.
\eeq
Here to ensure that $ \tilde L \in M_{6\times 2} (\bs )$
and $L \in M_{2\times 6} (\bs )$ we take $ u_+, v_+ , x_+ , y_+ \in \as_1$
and $ u_-, v_- , x_- , y_- \in \as_{-1}$.
Using the super determinant relation $ad - bc +  \hl_+ \hl_- = 1$,
one can verify that
\beq
u_+ = d, ~~~v_+ =-b,
~~~x_+ = (1+ 3\hl_+ \hl_- )b, ~~~y_+ =  (1+ 3\hl_+ \hl_- )d,
\eeq
\beq
u_- = -c, ~~~v_-= a,
~~~x_- =(1+ \hl_+ \hl_- )a, ~~~y_- = (1+ \hl_+ \hl_- )c,
\eeq
is a solution of (\ref{ll}), as needed.
\epf
By analogy with the classical situation,
we call $\as_{-1}$ and $\as_{1}$ the super-spin-bundle modules.
Proposition~\ref{free} is a super version of the fact that
the module of Dirac spinors,
i.e., the direct sum of the spin-bundle modules,
is free both for the classical and quantum sphere \cite{lps}.
In fact, the freeness of the module $P_{-1} \oplus P_1$ \cite[p.257]{hm99}
of Dirac spinors on the quantum sphere
can be shown by precisely the same method as in the super-sphere case.
It suffices to take in the proof of Proposition~\ref{free}
$M = \mbox{\scriptsize$\pmatrix{0 & 1 \cr 1 & 0 }$}$,
$F_{-1} = (\ha, \hg )^T (\hd, -q\hb )$,
$F_{1} = (\hd, \hb )^T (\ha, -q^{-1}\hg )$,
\beq
f_- = \pmatrix{\ha\cr \hg}  \pmatrix{\hd & -q\hb}~, ~~
f_+ = \pmatrix{\hb\cr \hd}  \pmatrix{-q^{-1}\hg & \ha}~,
\eeq
\beq
g_- = \pmatrix{\ha\cr \hg}  \pmatrix{\hd& -q\hb} ~, ~~
g_+ = \pmatrix{\hb\cr \hd}  \pmatrix{-q^{-1}\hg& \ha}~,
\eeq
where $\ha, \hb, \hg, \hd$ are the generators of \slq\ as in \cite{hm99}.

\section{Appendix: Gauge transformations}
\setcounter{equation}{0}

We follow here the definition of a gauge transformation used in
\cite{h-pm96}. For an $H$\gal\ $B\inc P$, it is defined as a unital
convolution-invertible homomorphism $f : H\ra P$ satisfying
$\dr\ci f = (f\ot \id )\ci \ad$.
We treat  this definition as the first approximation
of an appropriate concept of gauge transformations on \hge s
(see \cite{d-mb} and the paragraph above Proposition~3.4 in
\cite{h-pm96}, cf. \cite[Sections~6.1-6.2]{d-m97a},
\cite{d-m97b}, \cite[Section~5]{b-t96}).
It turns out that the space of strong connections is closed
under the action of gauge transformations \cite[Proposition~3.7]{h-pm96}.
The following theorem describes this action.
\bth\label{claim2}
Let $B\inc P$ be an $H$\gal ~admitting a strong connection.
The following describes a left action of gauge transformations
on strong connections which is compatible with the identifications
of Theorem~\ref{claim1}:
\vspace*{2.5mm}\\
1)
$(f\vt s) (p):= s \llp p\0 f(p\1 )\lrp f^{-1} (p\2 )$
\vspace*{2.5mm}\\
2) $(f\vt D) (p):= D \llp p\0 f(p\1 )\lrp f^{-1} (p\2 )$
\vspace*{2.5mm}\\
3) $(f\vt \Pi) (r\d p):= r \Pi \llp\d (p\0 f(p\1 ))\lrp f^{-1} (p\2 )
+ r p\0 f(p\1 )\d f^{-1} (p\2 ) $
\vspace*{2.5mm}\\
4) $(f\vt \ho ) (h):= f(h\1 )\ho (h\2 ) f^{-1} (h\3 )
+ f(h\1 )\d f^{-1} (h\2 ) $
\ethe
\bpf
We need to study the following diagrams:
\beq\label{diags}
\def\normalbaselines{\baselineskip30pt
\lineskip3pt \lineskiplimit3pt }
\def\mapright#1{\smash{
\mathop{\!\!\!-\!\!\!\longrightarrow\!\!\!}
\limits^{#1}}}
\def\mapdown#1{\Big\downarrow
\rlap{$\vcenter{\hbox{$\scriptstyle#1$}}$}}
\matrix{
GT(P)\times V_i & \mapright{\alpha_i}   & V_i \cr
\mapdown{{\rm {id}}\times J_{ij}} && \mapdown{J_{ij}}  \cr
GT(P)\times V_j & \mapright{\alpha_j}   & V_j \cr
}.
\eeq
Here $\ha_i$'s are the corresponding left actions specified above
and $J_{ij}$'s, $i,j \in \{1,2,3,4\}$ are obtained in an obvious
way by composing  suitable bijections $J_i$
introduced in the proof of Theorem~\ref{claim1}.
We know that $\ha_4$ is a well-defined left action
\cite[Proposition~3.4]{h-pm96}.
It suffices to show that
\beq
\ha_i = J_{4i}\ci \ha_4 \ci (\id\times J_{4i} )~~{\rm for}~~i \in \{1,2,3,\} .
\eeq
For $ i=3$ it is proved in \cite[Proposition~3.5]{h-pm96}.
For  $ i=2$, we have
\bea
&&\phantom{=}
J_{42} \llp f\vt J_{24}(D) \lrp (p)
\nonumber\\ &&
=(J_{12}\ci J_{41} )\llp f\vt (J_{34}\ci J_{23})(D) \lrp (p)
\nonumber\\ &&
=1\ot p  - J_{41}\llp f\vt (J_{34}\ci J_{23})(D) \lrp (p)
\nonumber\\ &&
=1\ot p - p\ot 1 - p\0 \llp f\vt (J_{34}\ci J_{23})(D) \lrp (p\1 )
\nonumber\\ &&
=\d p - p\0 f(p\1 )\llp J_{34}\ci J_{23} \lrp (D) (p\2 )f^{-1}(p\3 )
- p\0 f(p\1 )\d f^{-1}(p\2 )
\nonumber\\ &&
=\d  \llp p\0 f(p\1 )\lrp f^{-1}(p\2 )
- p\0 f(p\1 ) {p\2}\o J_{23}(D) (\d {p\2}\t )f^{-1}(p\3 ) \ .
\eea
Note now that, since $P$ admits a strong connection,
it is projective (Corollary~\ref{ff}) and hence flat as a left $B$-module.
Consequently $P\ot H$ is left $B$-flat
and
\beq
\ker \llp (\dr - \id \ot 1) \ot\sb B \id\ot\id \lrp = B\ot\sb B P\ot H\ .
\eeq
Using property (\ref{t1b}) of the translation map
and the $\ad$-colinearity of $f$, we obtain
\beq
(\dr - \id \ot 1)
\llp p\0 f(p\1 ) {p\2}\o \lrp \ot_B {p\2}\t \ot p\3 = 0 \ .
\eeq
Hence
\beq\label{flat}
p\0 f(p\1 ) {p\2}\o \ot_B {p\2}\t \ot p\3 \in B\ot\sb B P\ot H ,
\eeq
and we have:
\bea
J_{42} \llp f\vt J_{24}(D) \lrp (p)
\!\!\!\!\!\!&&
=\d  \llp p\0 f(p\1 )\lrp f^{-1}(p\2 )
-(J_{23}(D)\ci\d ) \llp p\0 f(p\1 ) {p\2}\o {p\2}\t \lrp f^{-1}(p\3 )
\nonumber\\ &&
=\d  \llp p\0 f(p\1 )\lrp f^{-1}(p\2 )
+(D-\d ) \llp p\0 f(p\1 ) \lrp f^{-1}(p\2 )
\nonumber\\ &&
=D \llp p\0 f(p\1 )\lrp f^{-1} (p\2 )
\nonumber\\ &&
= \ha_2(f, D) (p) \ .
\eea
Similarly, we compute:
\bea
J_{41} \llp f\vt J_{14}(s) \lrp (p)
\!\!\!\!\!\!&&
=p\ot 1 + p\0 f(p\1 )J_{14}(s) (p\2 ) f^{-1}(p\3 )
+ p\0 f(p\1 ) \d f^{-1}(p\2 ).
\eea
On the other hand,
\bea
J_{14}(s) (h)
\!\!\!\!\!\!&&
=  (J_{34}\ci J_{23}\ci J_{12})(s) (h)
\nonumber\\ &&
= h\o (J_{23}\ci J_{12})(\d h\t )
\nonumber\\ &&
= h\o (\d -J_{12}(s)) (h\t )
\nonumber\\ &&
= h\o (s-\id\ot 1) (h\t )
\nonumber\\ &&
= h\o s (h\t ) -\eps (h)\ot 1 \ .
\eea
Therefore, taking advantage of the left $B$-linearity of $s$, (\ref{flat}) and
(\ref{t3}), we obtain
\bea
&&
~~~
J_{41} \llp f\vt J_{14}(s) \lrp (p)
\nonumber\\ &&
=p\0 f(p\1 )\ot f^{-1}(p\2 )
+ p\0 f(p\1 ) {p\2}\o s({p\2}\t)f^{-1}(p\3 )
-p\0 f(p\1 )\ot f^{-1}(p\2 )
\nonumber\\ &&
= s \llp p\0 f(p\1 )\ot f^{-1}(p\2 )\lrp
\nonumber\\ &&
= \ha_1(f, s) (p) \ ,
\eea
as needed.
\epf
\vspace{-1cm}
\bre\label{pm}\em
The gauge transformations on $H$\gal\ $B\inc P$ are in on-to-one
correspondence with the gauge automorphisms
understood as unital left $B$-linear right $H$-colinear automorphisms of $P$
\cite[Proposition~5.2]{b-t96}.
If $f : H\ra P$ is a gauge  transformation,
then $F:P\ra P$, $F(p):= p\0 f(p\1 )$ is a  gauge automorphism.
Analogously, for $\ha\in \hO^1\! P$, we put
$F(\ha ) := \llp (\id\ot m)\ci (\id\ot\id\ot f )\ci \dsr \lrp (\ha )$.
(The other way round we have $f(h) = h\o F(h\t )$.)
Due to the right $H$-colinearity of the covariant differential $D$,
we can re-write point 2) of the above theorem as
$(D \triangleleft F ) (p) = F^{-1} (D F(p))$
This formula coincides with the usual formula for the action of
gauge transformations on  projective-module connections 
(e.g., see \cite[p.554]{c-a94}).
\ere
\bre\em
In the sense of the definition considered here,
the connections in Example~\ref{q} are not gauge equivalent.
This is because, for the quantum Hopf fibration,
any gauge transformation $f$ acts trivially on the space  of connections.
Indeed, since $H$ is spanned by group-like
elements, $f$ is convolution-invertible,
and the only invertible elements in \slq\ are non-zero complex numbers
\cite[Appendix]{hm99}, $f$ must be $\IC \setminus \{0\}$-valued.
This effect is due to working with non-completed (polynomial) algebras.
\ere

\footnotesize

\noindent {\it Acknowledgements.}
P.\ M.\ H.\ was partially supported
by  the  CNR postdoctoral fellowship at SISSA, Trieste, and KBN grant 2 P03A 030 14.
This work is part of H.\ G.'s project P11783-PHY
of the ``Fonds zur F{\"o}rderung der
Wissenschaftlichen Forschung in {\"O}sterreich."
H.\ G.\ thanks SISSA for an invitation to Trieste.
It is a pleasure to thank T.\ Brzezi\'nski, G.\ Landi, J.-L.\ Loday,
P.\ Schauenburg and A.\ Sitarz for very helpful discussions and suggestions.

\end{document}